\documentclass[a4paper,12pt]{elsarticle}
\usepackage{mathtools}
\usepackage{dsfont}
\usepackage{amsfonts}

\usepackage{mathrsfs}
\usepackage{geometry}
\usepackage{xcolor}
\definecolor{red}{RGB}{255,0,0} 
\definecolor{blue}{rgb}{0.0, 0.4, 0.65}
\definecolor{orange}{RGB}{255,84,0}
\definecolor{green}{RGB}{154,205,50}

\usepackage[overload]{empheq}

\usepackage{dimnum}
\usepackage{subfig}
\usepackage{hyperref}
\usepackage[capitalise]{cleveref}

\usepackage{multirow}

\usepackage{arydshln}

\newcommand{\B}[1]{\boldsymbol{#1}}
\newcommand{\FIN}{\textrm{fin}}
\newcommand{\LL}{\textsc{l}}
\newcommand{\LA}{\textsc{la}}
\renewcommand{\AA}{\textsc{a}}
\newcommand{\SF}{\textsc{sf}}
\newcommand{\SL}{\textsc{sl}}
\newcommand{\SA}{\textsc{sa}}

\newcommand{\W}{\textsc{w}}
\newcommand{\EQ}[1]{\eqref{eq:#1}}
\newcommand{\FIG}[1]{\ref{fig:#1}}
\newcommand{\TAB}[1]{\ref{tab:#1}}
\newcommand{\IR}{\mathbb{R}}
\newcommand{\dd}{\mathrm{d}}

\begin{document}

\begin{frontmatter}
\title{On the Consistency of Dynamic Wetting Boundary Conditions for the Navier--Stokes--Cahn--Hilliard Equations}

\author[TUE1]{T.H.B. Demont}
\author[TUE1]{S.K.F. Stoter}
\author[UT]{C. Diddens}
\author[TUE1]{E.H.~van~Brummelen}

\address[TUE1]{%
  Eindhoven University of Technology, Department of Mechanical Engineering,
  P.O.\ Box 513, 5600 MB Eindhoven, The Netherlands}
\address[UT]{
  University of Twente, Faculty of Science \& Technology, P.O.\ Box 217, 7500 AE Enschede, The Netherlands}

\begin{abstract}
We investigate the limiting behavior of the Navier--Stokes--Cahn--Hilliard model for binary-fluid flows as the diffuse-interface thickness passes to zero, in the presence of fluid--fluid--solid contact lines. Allowing for motion of such contact lines relative to the solid substrate is required to adequately model multi-phase and multi-species fluid transport past and through solid media. Even though diffuse-interface models provide an inherent slip mechanism through the mobility-induced diffusion, this slip vanishes as the interface thickness and mobility parameter tend to zero in the so-called sharp-interface limit. The objective of this work is to present dynamic wetting and generalized Navier boundary conditions for diffuse-interface models that are consistent in the sharp-interface limit. We concentrate our analysis on the prototypical binary-fluid Couette-flow problems. To verify the consistency of the diffuse-interface model in the limit of vanishing interface thickness, we provide reference limit solutions of a corresponding sharp-interface model. For parameter values both at and away from the critical viscosity ratio, we present and compare the results of both the diffuse- and sharp-interface models. The close match between both model results indicates that the considered test case lends itself well as a benchmark for further research.
\end{abstract}

\begin{keyword}
Navier--Stokes--Cahn--Hilliard equations\sep
sharp-interface limit\sep
dynamic wetting boundary condition\sep
generalized Navier boundary condition.
\end{keyword}

\end{frontmatter}

\section{Introduction}
\label{sec:intro}

Multi-phase fluid flows in which the individual fluid components are segregated by a small transition layer of molecular scale are ubiquitous in numerous fields of science and engineering. Inkjet printing is one such example of fields of application~\cite{lohse2022fundamental}. Advanced numerical models and simulation techniques, alongside experimental investigations, are used to accurately predict and simulate the complex physical phenomena that arise within these fields, providing valuable insight into the behavior of systems, and help the optimization of processes. There are generally two categories that mathematical-physical models for multi-phase fluid flows are classified under: sharp-interface models, in which an explicit representation by a manifold of co-dimension one separates the fluid components; and diffuse-interface models, where a thin-but-finite transition layer represents the interface between the fluid components. Sharp-interface models are of free-boundary type: the manifold accommodates kinematic and dynamic interface conditions which govern the advancement of the manifold by acting as coupled boundary conditions between adjacent fluid components' initial boundary-value problems. In diffuse-interface models, the transition layer across the pure fluid components consists of a mixture with a fluid-proportion that varies continuously and monotonically. Contrary to their sharp-interface counterpart, diffuse-interface models implicitly account for topological changes in the fluid-fluid interfaces, such as is the case for droplet coalescence or break-up~\cite{Lowengrub:1998uq}, and can intrinsically account for wetting scenarios where the fluid-fluid front propagates along a possibly elastic solid substrate~\cite{Yue:2010hq, Seppecher:1996kx, Jacqmin:2000kx, Brummelen:2021aw}.

Diffuse-interface models for two immiscible incompressible fluid species are generally described by the Navier--Stokes--Cahn--Hilliard (NSCH) equations. The NSCH equations represent a class of models rather than a single model, of which assorted variations have been presented over the past roughly 50 years. In the late 1970s, Hohenberg and Halperin presented the so-called ``model-H''~\cite{Hohenberg:1977hh}, where the densities of both fluid species are implied equal to each other. Approximately 20 years later, Lowengrub and Truskinovsky proposed a NSCH model~\cite{Lowengrub:1998uq} that allows for non-matching densities at the expense of being quasi-incompressible. Shokrpour et al.\ in~\cite{Simsek:2018gb} introduced a new quasi-incompressible NSCH system where the volume fraction takes on the role as order parameter. In the early 2010s, Abels, Garcke and Gr\"{u}n presented a thermodynamically consistent NSCH model (AGG NSCH), based on volume-averaged mixture velocity, for binary-fluids with non-matching densities~\cite{Abels:2012vn}. The AGG NSCH model, due to the premise of a volume-averaged mixture-velocity, has a solenoidal mixture velocity. Based on this property, on the fact of its thermodynamic consistency, and its consistent reduction to the underlying single-fluid Navier–Stokes equations, we turn to the AGG NSCH model in this paper, hereinafter referred to as the NSCH model.

Apart from the parameters related to the physical properties of the individual fluid components, NSCH models involve an additional three parameters related to the diffuse interface. Firstly, the interface-thickness parameter, $\varepsilon$, represents the orthogonal length scale of the transition layer between the two fluid components. Secondly, the mobility parameter, $m$, governs both the phase diffusion rate in the vicinity of the diffuse interface, and the rate at which the interface equilibrates in the phase-separated regime, thereby controlling the dynamics of the Ostwald-ripening effect. Introducing a phase-field dependence of the mobility parameter allows for further fine-tuning of the model behavior, but consequently introduces further complications in the numerical-approximation procedures~\cite{Barrett:1999nx}. Finally, the NSCH models contain a surface-tension parameter, $\sigma$, which is responsible for the excess free energy $\sigma_\LA$ of the diffuse interface according to $2\sqrt{2} \sigma = 3 \sigma_\LA$. While for the Navier--Stokes--Korteweg equations this proportionality depends on the interface thickness $\varepsilon$, this is not the case for the NSCH equations, where $\sigma$ is fully independent of $\varepsilon$ and $m$.


\subsection{Background}

When both $\varepsilon$ and $m$ tend to zero in the so-called sharp-interface limit $\varepsilon, m \rightarrow +0$, the transition layer collapses, under certain conditions, onto a manifold of co-dimension one. Even omitting situations involving moving contact lines and topological changes, the current understanding of the sharp-interface limit of the NSCH equations is limited. \cite[\S{}4.3]{Abels:2018ly} provides an overview of known results and open questions. Many authors have looked into the question regarding the appropriate scaling of the mobility parameter in relation to the interface-thickness parameter in the sharp-interface limit, which depends on the scaling $m \coloneqq m_\varepsilon$. If $m_\varepsilon$ vanishes suitably as $\varepsilon \rightarrow +0$, the classical sharp-interface Navier--Stokes binary-fluid model is obtained~\cite{Lowengrub:1998uq, Yue:2010hq, Abels:2018ly}. Contrarily, if $m_\varepsilon \propto \varepsilon^0$ as $\varepsilon \rightarrow +0$, then the NSCH model converges to the non-classical sharp-interface Navier--Stokes--Mullins--Sekerka model, see~\cite[\S{}4.1]{Abels:2018ly} and~\cite[\S{}4]{Jacqmin:2000kx}. Furthermore, if $m_\varepsilon = o(\varepsilon^3)$ as $\varepsilon \rightarrow +0$, then the resulting limit solution of the NSCH model generally violates the Young--Laplace condition on the pressure jump across the interface~\cite{Abels:2014ca}. Together, these results suggest that $m_\varepsilon$ satisfies a scaling with $\varepsilon$ proportional to $m_\varepsilon \propto \varepsilon^\alpha$ with $0 < \alpha \leq 3$ as $\varepsilon \rightarrow +0$, or some other more involved scaling relation satisfying these bounds. Several authors propose specific scaling relations or argue stricter bounds on $\alpha$~\cite{Khatavkar:2006gf, Jacqmin:2000kx}, with~\cite{Demont_Stoter_van_Brummelen_2023} suggesting an optimal scaling of $\alpha \approx 1.7$ for the case of oscillating droplet scenarios.

When introducing contact line dynamics to the problem at hand, further complications are induced that require resolving; specifically, when considering desired limiting behavior of the NSCH equations, additional care has to be taken in the imposing of adequate boundary conditions. In the early '70s, Huh and Scriven derive solution fields for the sharp interface case with no-slip boundary conditions~\cite{huh1971hydrodynamic}.  Using the Stokes equations, they describe the movement of the contact line over a flat solid surface and provide an excellent overview for different flow profiles at different viscosity ratios and contact angles, including so-called triple-wedge flows. Huh and Scriven conclude that while the velocity fields appear realistic, the stresses and viscous dissipation are unbounded when approaching the contact point. About 15 years later, in the mid '80s, Cox derives expressions relating the micro- and macroscopic contact angles using matched asymptotic expansions~\cite{Cox:1986jq}. He considers slip near the contact line according to a slip length $s_\nu$ in order to remove the stress singularity at the contact line. Another 15 years later, Jacqmin studies the contact line dynamics of the diffuse-interface NSCH equations with no-slip boundary conditions, yet with a boundary slip induced through diffusion governed by the mobility parameter $m$~\cite{Jacqmin:2000kx}. Jacqmin derives a corresponding diffusive slip length $s_m$ dependent on the mobility $m$ equivalent to Cox' $s_\nu$ according to $s_\nu \equiv s_m \coloneqq \sqrt{\eta m}$, where $\eta$ is the viscosity of both fluid components. In 2010, Yue et al.\ advocate a proportionality of $m$ to the interface thickness $\varepsilon$ according to $m \propto \varepsilon^0$~\cite{Yue:2010hq}. They show that convergence occurs as $\varepsilon \rightarrow +0$, where the results are independent of $\varepsilon$ when $\varepsilon < 4 s_m$. However, as emphasized in the preceding paragraph, one would desire an $m_\varepsilon = o(\varepsilon^0)$ scaling as $\varepsilon \rightarrow +0$, i.e.\ $m_\varepsilon \propto \varepsilon^\alpha$ with $\alpha > 0$, in order to obtain the classical NS sharp-interface equations, which is in direct contradiction with a constant diffusive slip-length, and thus constant $m$. In earlier work, we conjecture that the generalized Navier boundary conditions can be employed as a way out of this `Morton's fork' \cite{Stoter2023b}, a concept we study systematically in the current article. The same approach is taken in Rougier et al., for a Stokes flow model using a level-set function and generalized Navier boundary conditions (GNBC)~\cite{rougier2021slip}. They investigate the correct slip length $s_\nu$ by comparing to experiments, and show that in order to reproduce experimentally measured contact angles, $s_\nu$ needs to depend on the contact line velocity.


\subsection{Overview}

In this article, we address the issue of undesired convergence behavior of diffuse-interface models in the presence of contact line dynamics. Given the results on the $m_\varepsilon$ scaling in the sharp-interface limit, extra constraints on the scaling of $m_\varepsilon$ are undesired when including contact line dynamics. Furthermore, convergence in the sharp-interface limit to the classical sharp-interface solution is quintessential. As such, we introduce the generalized Navier boundary conditions, accompanied by dynamic contact angle boundary conditions, to the wetting boundaries of the NSCH model. We study the properties of the model and its sharp-interface limit through various test-cases.

In particular, we compare the results of the no-slip NSCH model with purely mobility-induced diffusive slip $s_m$ to a sharp-interface counterpart that incorporates the slip length $s_\nu$ in the generalized Navier boundary conditions, where the respective diffusive and GNBC slip lengths are equated according to \cite{Cox:1986jq, Jacqmin:2000kx}. We investigate the convergence as $\varepsilon \rightarrow +0$, while fixing the diffusive slip length $s_m$ and thus fixing the mobility $m$, where we consider $\ell \gg s_{\nu,m} \gg \varepsilon$, with $\ell$ the macroscopic characteristic length scale. Furthermore, we look into a range of different slip-lengths and capillary numbers. Likewise, we look into the sharp-interface limit $\varepsilon, m \rightarrow +0$ with no-slip boundary conditions and thus a vanishing diffusive slip length $s_m$, where we expect an ill-defined limit solution and a singularity in the shear force in a small window near the contact line \cite{huh1971hydrodynamic}. Next, we introduce a positive slip length $s_\nu$ through the generalized Navier boundary condition in the NSCH model in order to resolve the inconsistency in the limit and once again look at the convergence as $\varepsilon, m \rightarrow +0$, where $s_m \ll s_\nu$, to the sharp-interface limit with equal GNBC slip length $s_\nu$. By considering separate cases, where we either deviate from a 1:1 fluid viscosity ratio or deviate from a 90 degree microscopic contact angle, we induce the triple-wedge flows as observed in \cite{huh1971hydrodynamic}.

In order to conduct our analysis of the sharp-interface limit of the NSCH equations with generalized Navier boundary conditions in the setting of contact line dynamics, we consider two-dimensional Couette-type flow problems. Several interface-related quantities are measured and used as measurables with which to compare to the theorized sharp-interface limits, such as the contact line displacement and macroscopic contact angle. An adaptive finite-element method, in which the adaptivity is guided by an a posteriori error estimate (see~\cite{Brummelen:2021aw, Demont:2022dk} for details) is applied to enable an exploration of the asymptotic regime. With the close match between the results of the sharp-interface model and that of the limiting diffuse-interface model, we propose the considered test cases as a benchmark for further research: both the diffuse- and sharp-interface models have been implemented fully independent from each other.

The structure of the remainder of this paper is as follows. In Section~\ref{sec:prob_form}, we present the diffuse-interface Abels--Garcke--Gr\"{u}n NSCH model, along with the wetting and generalized Navier boundary conditions. In Section~\ref{sec:model_prop}, we show the model is thermodynamically consistent and perform an analysis on the dynamic contact angle. In Section~\ref{sec:num_exp}, we summarize the approximation setting, discretization schemes, and the various test cases, and conduct the numerical experiments in order to study the NSCH model with generalized Navier boundary conditions in the setting of contact line dynamics. Finally, in Section~\ref{sec:concl} we conclude with a summary and a discussion.


\section{Governing equations}
\label{sec:prob_form}

\begin{figure}[!b]
\centering
\includegraphics{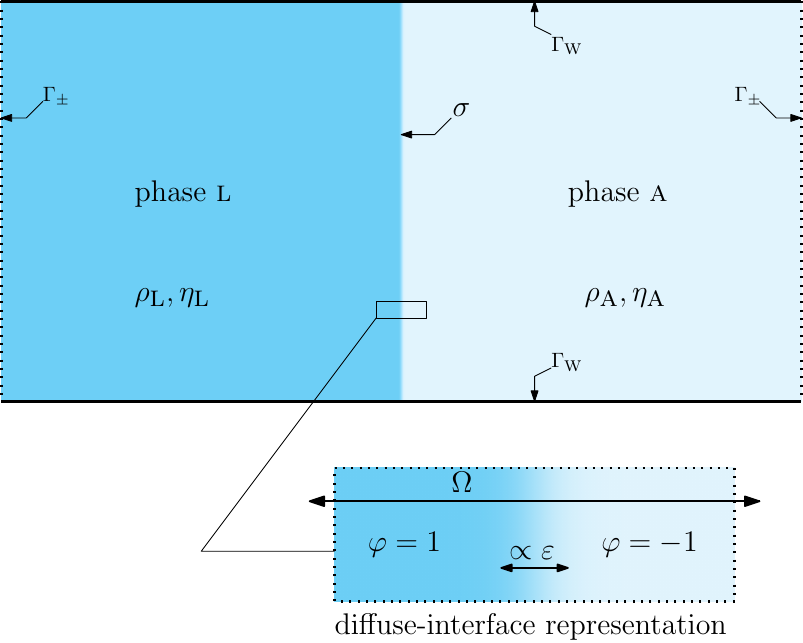}
\caption{A schematic of the considered physical setting of two fluid regions in a thin channel is displayed, located around the fluid-fluid interface. The zoom highlights the diffuse-interface representation of the setting, with the indicated phase-field variable and interface thickness.}
\label{fig:channel_model}
\end{figure}

We consider a binary fluid consisting of two incompressible, isothermal, immiscible, and Newtonian species with non-zero viscosity. As the (initial) setting for the binary fluid under consideration, we define two fluid regions in a thin channel, separated by a finite-thickness fluid-fluid interface comprised of a mixture of both fluids; see the illustration in Figure \ref{fig:channel_model}. In order to notationally adhere to the conventional setup of a liquid submerged in an ambient, we stick to the naming convention of labeling the two fluid domains as $\LL$ for liquid, and $\AA$ for ambient.


\subsection{Diffuse-interface representation}
We utilize the Abels--Garcke--Gr\"{u}n NSCH equations on an open time interval $(0, t_\FIN) \subseteq \mathds{R}_{>0}$ and a spatial domain corresponding to a simply connected time-independent subset $\Omega \subseteq \mathds{R}^d$ ($d=2,3$), in order to describe the evolution of the phase field parameter $\varphi \in [-1,1]$ representing pure species $\LL$ and $\AA$ when $\varphi=1$ and $\varphi=-1$, respectively, and a mixture of both when $\varphi \in (-1,1)$. Besides the order parameter $\varphi$, the chemical potential $\mu$, the volume-averaged velocity $\B{u}$, and the pressure $p$ are the unknown fields. The NSCH model's bulk equations, as presented in~\cite{Abels:2012vn}, are given by
\begin{subequations}
\label{eq:strong}
\begin{empheq}[right={\empheqrbrace \textrm{ in } \Omega}]{align}
\partial_t \left( \rho \B{u} \right) + \nabla \cdot \left( \rho \B{u} \otimes \B{u} \right) + \nabla \cdot \left( \B{u} \otimes \B{J} \right) + \nabla p - \nabla \cdot \B{\tau} - \nabla \cdot \B{\zeta} & = 0, \label{eq:strong_sub1}\\
\nabla \cdot \B{u} & = 0, \label{eq:strong_sub2}\\
\partial_t \varphi + \nabla \cdot \left( \varphi \B{u} \right) - \nabla \cdot \left( m \nabla \mu \right) & = 0, \label{eq:strong_sub3}\\
\mu + \sigma \varepsilon \Delta \varphi - \frac{ \sigma } { \varepsilon } \Psi' & = 0. \label{eq:strong_sub4}
\end{empheq}
\end{subequations}
The relative mass flux $\B{J}$, the viscous stress $\B{\tau}$, the capillary stress $\B{\zeta}$, and the mixture energy density $\Psi$ all require closure relations, which are defined as
\begin{subequations}
\label{eq:closureeqs}
\begin{alignat}{3}
& \B{J} \coloneqq m \frac{ \rho_{\AA} - \rho_{\LL} } { 2 } \nabla \mu  \,, \label{eq:closureeqs_sub1}\\
& \B{\tau} \coloneqq \eta ( \nabla \B{u} + ( \nabla \B{u} )^T ) \,, \label{eq:closureeqs_sub2}\\
& \B{\zeta} \coloneqq - \sigma \varepsilon \nabla \varphi \otimes \nabla \varphi + \B{I} \left( \frac{ \sigma \varepsilon } { 2 } | \nabla \varphi |^2 + \frac{ \sigma } { \varepsilon } \Psi \right) \,, \label{eq:closureeqs_sub3}\\
& \Psi \left( \varphi \right) \coloneqq \frac{1}{4} \left( \varphi^2 -1 \right)^2 \,. \label{eq:closureeqs_sub4}
\end{alignat}
\end{subequations}
The interface thickness parameter $\varepsilon>0$ and mobility parameter $m>0$ are the model parameters, which govern the length and time scale of the diffuse interface, respectively. Finally, the material parameters are $\sigma$ --- a rescaling of the droplet-ambient surface tension $\sigma_\LA$ according to $2\sqrt{2}\sigma=3\sigma_\LA$ --- and the mixture density $\rho$ and mixture viscosity $\eta$, which generally depend on~$\varphi$. The phase field parameter $\varphi$ must be allowed to take on values outside of $[-1,1]$ in order to ensure existence of a solution to the system of equations~\cite{Grun:2016gi}. To ensure positive densities even for the nonphysical scenario in which $\varphi \notin [-1,1]$, the mixture density function $\rho(\varphi)$ includes an extension~\cite{Bonart:2019re}:
\begin{equation}\label{eq:densityextension}
\rho(\varphi) = \left\{
\begin{tabular}{ll}
$\frac{ 1 } { 4 } \rho_{\AA},$ & $\varphi \leq - 1 - 2 \lambda \,,$\\
$\frac{ 1 } { 4 } \rho_{\AA} + \frac{ 1 } { 4 } \rho_{\AA} \lambda^{-2} \left( 1 + 2 \lambda + \varphi \right)^2,$ & $\varphi \in ( - 1 - 2 \lambda , - 1 - \lambda) \,,$\\
$\frac{ 1 + \varphi } { 2 } \rho_{\LL} + \frac{ 1 - \varphi } { 2 } \rho_{\AA},$ & $\varphi \in [ - 1 - \lambda, 1 + \lambda ] \,,$\\
$\rho_{\LL} + \frac{ 3 } { 4 } \rho_{\AA} - \frac{ 1 } { 4 } \rho_{\AA} \lambda^{-2} \left( 1 + 2 \lambda - \varphi \right)^2,$ & $\varphi \in ( 1 + \lambda, 1 + 2 \lambda ) \,,$\\
$\rho_{\LL} + \frac{ 3 } { 4 } \rho_{\AA},$ & $\varphi \geq 1 + 2 \lambda \,,$
\end{tabular}
\right.
\end{equation}
where $\rho_\LL$ and $\rho_\AA$ --- which are present in the definition of the relative mass flux $\B{J}$ as well --- are the liquid and ambient fluid densities, respectively. $\lambda \coloneqq \rho_{\AA} / \left( \rho_{\LL} - \rho_{\AA} \right)$ is the relative threshold parameter beyond which point the positivity-preserving density extention takes effect. Additionally, we apply the Arrhenius mixture-viscosity model~\cite{Arrhenius:1887xr} for the viscosity interpolation:
\begin{equation}\label{eq:arrhenius}
\log \eta( \varphi ) = \frac{ \left( 1 + \varphi \right) \Lambda \log \eta_{\LL} + \left( 1 - \varphi \right) \log \eta_{\AA} } { \left( 1 + \varphi \right) \Lambda + \left( 1 - \varphi \right) } \,,
\end{equation}
where $\Lambda \coloneqq \frac{ \rho_{\LL} M_{\AA} } { \rho_{\AA} M_{\LL} }$ is the intrinsic volume ratio between the two fluids, and $M_{\LL}$ and $M_{\AA}$ are their respective molar masses. In order to eliminate singularities, we choose an intrinsic volume ratio of $\Lambda=1$ in this work; see \emph{Remark 2} in \cite{Brummelen:2021aw} for further details.


\subsection{Wetting boundary conditions}
\label{sec:wetting_BC}

The considered physical setting of laminar two-phase flow through a small channel allows for an adequate investigation of the contact line dynamics, while minimizing the complexity --- for small enough \Ca[l] and \Rey[l] --- of the dynamics in the bulk. We then enforce on the wetting boundaries $\Gamma_\W$ of the NSCH model the dynamic contact angle boundary conditions and the generalized Navier boundary conditions~\cite{Stoter2023b}, according to
\begin{subequations}
\label{eq:wetting_BC}
\begin{empheq}[right={\empheqrbrace \textrm{ on } \Gamma_\W}]{align}
\partial_t \varphi + \B{u}_S \cdot \nabla_S \varphi & = -\nu_1 \left( \sigma \varepsilon \partial_n \varphi + \sigma'_{\SF}(\varphi) \right)\,, \label{eq:wetting_BC_sub1}\\
\B{u}_S - \B{u}_\W & = \nu_2 \left( (p\B{n} - \B{\tau}\B{n} - \B{\zeta}\B{n})_S + \nabla_S \sigma_\SF(\varphi) \right)\,, \label{eq:wetting_BC_sub2}\\
u_n & = 0\,, \label{eq:wetting_BC_sub3}\\
\partial_n \mu & = 0\,, \label{eq:wetting_BC_sub4}
\end{empheq}
\end{subequations}
where $\B{n}$ denotes the exterior unit normal vector on~$\partial\Omega$, $\B{u}_\W$ is the external boundary wall velocity, $(\cdot)_{S}=\B{n}\times(\cdot)\times\B{n}$ is the in-plane component of a vector, and $\sigma_\SF$ represents the solid-fluid surface tension according to
\begin{equation}
\sigma_\SF(\varphi) = \frac{1}{4} \left( \varphi^3 - 3 \varphi \right) \left( \sigma_\SA - \sigma_\SL \right) + \frac{1}{2} \left( \sigma_\SL + \sigma_\SA \right),
\end{equation}
with $\sigma_\SL \geq 0$ and $\sigma_\SA \geq 0$ the solid-fluid surface tensions of phase $\LL$ and $\AA$, respectively; see~\cite{Brummelen:2016qa, Jacqmin:2000kx, Shokrpour-Roudbari:2016dp, Yue:2011uq}. By coupling these wetting boundary conditions to the bulk equations, there is direct control over the (static) equilibrium microscopic contact angle that the fluid-fluid interface makes with the solid boundary via $\sigma_\SL$ and $\sigma_\SA$. Furthermore, there is a freedom in choice between no, limited, and free slip --- aside from inherent slip through the mobility parameter --- local to the wetting boundary via the $\nu_2$ parameter ($\nu_2 = 0$, $\nu_2 > 0$, and $\nu_2~\rightarrow~\infty$, respectively). Finally, additional control over the microscopic contact angle behavior is introduced via the $\nu_1$ parameter.

Additionally, at the far ends of the channel domain, $\Gamma_\pm$, we will impose in- and outflow boundary conditions according to
\begin{subequations}
\label{eq:pm_BC}
\begin{empheq}[right={\empheqrbrace \textrm{ on } \Gamma_\pm}]{align}
\B{u} & = \overline{\B{u}}_{\pm}\,,\\
\partial_n\varphi&=0\,,\\
\partial_n\mu&=0\,,
\end{empheq}
\end{subequations}
with Neumann boundary conditions for the order parameter $\varphi$ and chemical potential $\mu$. The in- and outflow profile $\overline{\B{u}}_{\pm}$ will be further specified in Section \ref{sec:setup}, where the exact details of the considered test cases are outlined. Note that the additional constraint of a homogeneous average pressure is imposed for consistency.


\section{Properties}
\label{sec:model_prop}

In this section, several properties of the NSCH model \EQ{strong} coupled with the wetting equations~\EQ{wetting_BC} are presented. In Section~\ref{sec:nrg_dissip}, the energy dissipation structure of the coupled equations is established. Section~\ref{sec:SI_lim_analysis} presents a partial, formal sharp-interface limit analysis of the coupled equations. Specifically, we show that the Ginzburg--Landau energy density converges to a delta distribution on $\Gamma_\LA$, multiplied by liquid-ambient surface tension, as $\varepsilon \rightarrow +0$, and we establish that the diffuse-interface GNBC~\eqref{eq:wetting_BC_sub2} formally converges to the sharp-interface GNBC. Finally, in Section~\ref{sec:CA_condition}, we focus on the equilibrium contact-angle condition that emerges in the limit $\nu_1 \rightarrow \infty$. We show that if the microscopic contact angle is equal to the equilibrium contact angle, the phase field in the vicinity of the contact line is compatible with the equilibrium contact-angle condition, independent of $\varepsilon$. Additionally, if the contact-line velocity coincides with the flow velocity orthogonal to the contact line, the phase field in the vicinity of the contact line is compatible with the dynamic contact angle boundary condition~\eqref{eq:wetting_BC_sub1}, independent of~$\varepsilon$ as well.


\subsection{Energy dissipation}
\label{sec:nrg_dissip}

The energy functional $\mathscr{E}$ of the NSCH model~\EQ{strong} coupled with the wetting equations~\EQ{wetting_BC} consists of the kinetic energy of the fluid, the free energy associated with (mixing of) both fluid species, and the solid--fluid surface energy, and is given by
\begin{equation}
\label{eq:nrg_functnl}
\mathscr{E}(t) = \int_\Omega\left( \frac{\sigma\varepsilon}{2} | \nabla \varphi |^2 + \frac{\sigma}{\varepsilon} \Psi \right) + \int_\Omega \frac{1}{2}\rho |\B{u}|^2 + \int_{\Gamma_\W} \sigma_\SF(\varphi)\,.
\end{equation}
Here, the first term is a Ginzburg--Landau-type energy functional associated with the (mixing of) the binary-fluid components, the second term corresponds to the kinetic energy of the binary fluid, and the final term represents the energetic contribution of the solid--fluid surface tension. In the following, we show that solutions that satisfy the equations as presented in Section~\ref{sec:prob_form}, adhere to the thermodynamic consistency condition:
\begin{equation}
\label{eq:nrg_dissip}
\begin{aligned}
\frac{d}{dt} \mathscr{E}(t)  =& - \int_\Omega m |\nabla \mu|^2 - \int_\Omega \nabla \B{u} : \B{\tau} - \int_{\Gamma_\W} \nu_1 \left| \sigma \varepsilon \partial_n \varphi + \sigma'_\SF(\varphi) \right|^2
                            \\
                             & 
                             - \int_{\Gamma_\W} \nu_2^{-1} | \B{u}_S - \B{u}_\W |^2+ \int_{\Gamma_\W} \nu_2^{-1}\B{u}_\W \cdot (\B{u}_\W-\B{u}_S)+ \textrm{bnd}_\pm\\
                            \leq & \int_{\Gamma_\W} \nu_2^{-1}\B{u}_\W  \cdot (\B{u}_\W-\B{u}_S) + \textrm{bnd}_\pm\,,
\end{aligned}
\end{equation}
where `bnd$_\pm$' is a placeholder term encompassing terms that only have support on the in- and outflow boundaries $\Gamma_\pm$ of the binary-fluid system. This implies that the solutions have a free-energy dissipation property, and the only potential localized energy inputs are through the in- and outflow boundaries $\Gamma_\pm$ and due to a non-zero wetting boundary wall velocity $\B{u}_\W$. It is noteworthy that the third and fourth terms after the identity are dissipative terms emerging from the wetting boundary conditions~\eqref{eq:wetting_BC_sub1}--\eqref{eq:wetting_BC_sub2}, while the fifth term is an energy-production term emerging from the generalized Navier boundary condition~\eqref{eq:wetting_BC_sub2}, associated with the work exerted on the binary fluid by the moving wall. Below, we will consider each individual term in the energy~\EQ{nrg_functnl} to prove the energy dissipation relation~\EQ{nrg_dissip}.

For the Ginzburg--Landau-type energy functional, we obtain
\begin{equation}
\label{eq:GL_int1}
\begin{aligned}
& \frac{d}{dt} \int_\Omega\left( \frac{\sigma\varepsilon}{2} | \nabla \varphi |^2 + \frac{\sigma}{\varepsilon} \Psi \right)\\
& \quad = \int_\Omega \sigma \varepsilon \nabla \varphi \cdot \nabla \partial_t \varphi + \frac{\sigma}{\varepsilon} \Psi' \partial_t \varphi\\
& \quad = \int_\Omega \partial_t \varphi \left( \frac{\sigma}{\varepsilon} \Psi' - \sigma \varepsilon \Delta \varphi \right) + \int_{\partial \Omega} \sigma \varepsilon \partial_t \varphi \partial_n \varphi.
\end{aligned}
\end{equation}
Substituting~\EQ{strong_sub3} and~\EQ{strong_sub4}, using the product rule twice, and applying the divergence theorem, we find that
\begin{equation}
\label{eq:GL_int2}
\begin{aligned}
& \int_\Omega \partial_t \varphi \left( \frac{\sigma}{\varepsilon} \Psi' - \sigma \varepsilon \Delta \varphi \right)\\
& \quad = \int_\Omega \left( \nabla \cdot \left( m \nabla \mu \right) \right) \mu - \int_\Omega \mu \nabla \cdot \left( \varphi \B{u} \right)\\
& \quad = - \int_\Omega m |\nabla \mu|^2 + \int_{\partial \Omega} m \mu \partial_n \mu - \int_\Omega \B{u} \cdot \mu \nabla \varphi + \int_\Omega \mu \varphi \nabla \cdot \B{u}\,.
\end{aligned}
\end{equation}
Standard tensor-calculus operations along with a substitution of~\EQ{strong_sub4} yields the identities
\begin{equation}
\begin{aligned}
\nabla \cdot \left( \sigma \varepsilon \nabla \varphi \otimes \nabla \varphi \right) & = \sigma \varepsilon \left( \Delta \varphi \nabla \varphi + \frac{1}{2} \nabla | \nabla \varphi |^2 \right)\\
& = \nabla \varphi \left( \frac{\sigma}{\varepsilon} \Psi' - \mu \right) + \nabla \left( \frac{\sigma \varepsilon}{2} | \nabla \varphi |^2 \right)\\
& = \nabla \left( \frac{\sigma}{\varepsilon} \Psi + \frac{\sigma \varepsilon}{2} |\nabla \varphi|^2 \right) - \mu \nabla \varphi\,.
\end{aligned}
\end{equation}
Using these identities and definition~\EQ{closureeqs_sub3} for $\B{\zeta}$, along with the product rule and divergence theorem twice more, we find that
\begin{equation}
\begin{aligned}
- \int_\Omega \B{u} \cdot \mu \nabla \varphi & = - \int_\Omega \B{u} \cdot \left\{ \nabla \left( \frac{\sigma}{\varepsilon} \Psi + \frac{\sigma \varepsilon}{2} | \nabla \varphi |^2 \right) - \nabla \cdot \left( \sigma \varepsilon \nabla \varphi \otimes \nabla \varphi \right) \right\}\\
& = - \int_{\partial \Omega} u_n \left( \frac{\sigma}{\varepsilon} \Psi + \frac{\sigma \varepsilon}{2} |\nabla \varphi|^2 \right) + \int_\Omega \nabla \B{u} : \B{I} \left( \frac{\sigma}{\varepsilon} \Psi + \frac{\sigma \varepsilon}{2} |\nabla \varphi|^2 \right)\\
& \quad - \int_\Omega \nabla \B{u} : \left( \sigma \varepsilon \nabla \varphi \otimes \nabla \varphi \right) + \int_{\partial \Omega} \sigma \varepsilon \partial_n \varphi \B{u} \cdot \nabla \varphi\\
& = - \int_{\partial \Omega} u_n \left( \frac{\sigma}{\varepsilon} \Psi + \frac{\sigma \varepsilon}{2} |\nabla \varphi|^2 \right)\\
& \quad + \int_\Omega \nabla \B{u} : \B{\zeta} + \int_{\partial \Omega} \sigma \varepsilon \partial_n \varphi \B{u} \cdot \nabla \varphi\,.
\end{aligned}
\end{equation}
Substituting this result into~\EQ{GL_int2} and, in turn, into~\EQ{GL_int1}, and recalling the solenoidality property~\EQ{strong_sub2}, we ultimately find:
\begin{equation}
\label{eq:E_GL_res}
\begin{aligned}
\frac{d}{dt} \int_\Omega\left( \frac{\sigma\varepsilon}{2} | \nabla \varphi |^2 + \frac{\sigma}{\varepsilon} \Psi \right) = & - \int_\Omega m |\nabla \mu|^2 + \int_{\partial \Omega} m \mu \partial_n \mu\\
& - \int_{\partial \Omega} u_n \left( \frac{\sigma}{\varepsilon} \Psi + \frac{\sigma \varepsilon}{2} |\nabla \varphi|^2 \right) + \int_\Omega \nabla \B{u} : \B{\zeta}\\
& + \int_{\partial \Omega} \sigma \varepsilon \partial_n \varphi \B{u} \cdot \nabla \varphi + \int_{\partial \Omega} \sigma \varepsilon \partial_t \varphi \partial_n \varphi\,.
\end{aligned}
\end{equation}

Next, let us consider the kinetic-energy contribution. It holds that
\begin{equation}
\frac{\partial}{\partial t} \frac{1}{2} \rho |\B{u}|^2 + \nabla\cdot \B{u} \frac{1}{2} \rho | \B{u} |^2
= \B{u} \cdot \left( \partial_t \rho \B{u} + \nabla \cdot \left( \rho \B{u} \otimes \B{u} \right) \right)
- \frac{1}{2} | \B{u} |^2 \left( \partial_t \rho + \nabla \cdot \left( \rho \B{u} \right) \right)\,.
\end{equation}
The conservation-of-mass relation for the NSCH model is given by~\cite[\S{}3.1]{Abels:2018ly} 
\begin{equation}
\partial_t \rho + \nabla \cdot \left( \rho \B{u} \right) + \nabla \cdot \B{J} = 0\,,
\end{equation}
and together with the conservation-of-momentum equation~\EQ{strong_sub1} and the divergence theorem, we find that the kinetic-energy term can be recast as
\begin{equation}
\begin{aligned}
\frac{d}{dt} \int_\Omega \frac{1}{2}\rho |\B{u}|^2 & = \int_\Omega \B{u} \cdot \left( - \nabla \cdot \left( \B{u} \otimes \B{J} \right) - \nabla p + \nabla \cdot \B{\tau} + \nabla \cdot \B{\zeta} \right)\\
& \quad + \frac{1}{2} | \B{u} |^2 \nabla \cdot \B{J} - \int_{\partial \Omega} \frac{1}{2} \rho | \B{u} |^2 u_n\,.
\end{aligned}
\end{equation}
By applying tensor calculus operations, one can derive that
\begin{equation}
\label{eq:J_calc}
- \int_\Omega \B{u} \cdot \nabla \cdot \left( \B{u} \otimes \B{J} \right) + \int_\Omega \frac{1}{2} | \B{u} |^2 \nabla \cdot \B{J} = - \int_\Omega \nabla \cdot \left(\frac{1}{2} | \B{u} |^2 \B{J} \right)\,.
\end{equation}
After some rewriting, substitution of~\EQ{J_calc}, and applying the product rule and divergence theorem, we obtain for the kinetic-energy contribution the following expression:
\begin{equation}
\label{eq:E_kin_res}
\begin{aligned}
\frac{d}{dt} \int_\Omega \frac{1}{2}\rho |\B{u}|^2 =&  - \int_\Omega \nabla \B{u} : \B{\tau} - \int_\Omega \nabla \B{u} : \B{\zeta}
\\
 &+ \int_{\partial \Omega} \B{u} \cdot \left( \B{\tau} \B{n} + \B{\zeta} \B{n} - p \B{n} \right)
- \int_{\partial \Omega} \frac{1}{2} |\B{u}|^2 \left( \rho u_n + J_n \right)\,.
\end{aligned}
\end{equation}

Now, bringing the resulting Equations~\EQ{E_GL_res} and~\EQ{E_kin_res} together, we obtain
\begin{equation}
\label{eq:E_pre_final}
\begin{aligned}
\frac{d}{dt} \mathscr{E}(t) = & - \int_\Omega m |\nabla \mu|^2 + \int_{\partial \Omega} m \mu \partial_n \mu\\
& - \int_{\partial \Omega} u_n \left( \frac{\sigma}{\varepsilon} \Psi + \frac{\sigma \varepsilon}{2} |\nabla \varphi|^2 \right)\\
& + \int_{\partial \Omega} \sigma \varepsilon \partial_n \varphi \B{u} \cdot \nabla \varphi + \int_{\partial \Omega} \sigma \varepsilon \partial_t \varphi \partial_n \varphi\\
& - \int_\Omega \nabla \B{u} : \B{\tau} + \int_{\partial \Omega} \B{u} \cdot \left( \B{\tau} \B{n} + \B{\zeta} \B{n} - p \B{n} \right)\\
& - \int_{\partial \Omega} \frac{1}{2} |\B{u}|^2 \left( \rho u_n + J_n \right) + \int_{\Gamma_\W} \sigma'_\SF(\varphi) \partial_t \varphi\,.
\end{aligned}
\end{equation}
By adding the partition of zero
\begin{equation}
0 = \int_{\partial \Omega} \B{u} \cdot \nabla \sigma_\SF(\varphi) - \int_{\partial \Omega} \B{u} \cdot \nabla \sigma_\SF(\varphi)\,,
\end{equation}
to Equation~\EQ{E_pre_final}, and making use of the homogeneous boundary conditions in~\EQ{wetting_BC} and~\EQ{pm_BC},  we obtain --- after some rewriting and grouping terms that only have support on the in- and outflow boundaries $\Gamma_\pm$ under `bnd$_\pm$':
\begin{equation}
\begin{aligned}
\frac{d}{dt} \mathscr{E}(t) & = - \int_\Omega m |\nabla \mu|^2 - \int_\Omega \nabla \B{u} : \B{\tau} + \int_{\Gamma_\W} \B{u}_S \cdot \left( \B{\tau} \B{n} + \B{\zeta} \B{n} - p \B{n} \right)\\
& \quad - \int_{\Gamma_\W} \B{u}_S \cdot \nabla_S \sigma_\SF(\varphi) + \int_{\Gamma_\W} \B{u}_S \cdot \nabla_S \varphi \, \sigma \varepsilon \partial_n \varphi + \int_{\Gamma_\W} \partial_t \varphi \, \sigma \varepsilon \partial_n \varphi\\
& \quad + \int_{\Gamma_\W} \partial_t \varphi \, \sigma'_\SF(\varphi) + \int_{\Gamma_\W} \B{u}_S \cdot \nabla_S \varphi \, \sigma'_\SF(\varphi) + \textrm{bnd}_\pm\\
& = - \int_\Omega m |\nabla \mu|^2 - \int_\Omega \nabla \B{u} : \B{\tau}+ \textrm{bnd}_\pm\\
& \quad + \int_{\Gamma_\W} \left( \partial_t \varphi + \B{u}_S \cdot \nabla_S \varphi \right) \left( \sigma \varepsilon \partial_n \varphi + \sigma'_\SF(\varphi) \right)\\
& \quad + \int_{\Gamma_\W} \B{u}_S \cdot \left( \left( \B{\tau} \B{n} + \B{\zeta} \B{n} - p \B{n} \right) - \nabla_S \sigma_\SF(\varphi) \right)\,.
\end{aligned}
\end{equation}
Finally, substitution of the dynamic contact angle boundary condition~\EQ{wetting_BC_sub1} and the generalized Navier boundary condition~\EQ{wetting_BC_sub2} yields the desired energy dissipation property in~\EQ{nrg_dissip}, establishing thermodynamic consistency of~\eqref{eq:strong} with state functions~\eqref{eq:closureeqs} subject to boundary conditions~\eqref{eq:wetting_BC},with the only energy production being through the in- and outflow boundaries $\Gamma_\pm$ and by a non-zero velocity $\B{u}_\W$ of the wetting boundary wall. 


\subsection{Sharp-interface limit}
\label{sec:SI_lim_analysis}
In this section we present an analysis of the NSCH model~\eqref{eq:strong} with boundary conditions~\eqref{eq:wetting_BC} in the sharp-interface limit. Rigorous analyses of sharp-interface limits are technical, and we restrict ourselves here to formal derivations. Our presentation of the sharp-interface GNBC follows that of~\cite{Gerbeau:2009km}.

To elucidate our approach, we regard the stationary Cahn--Hilliard equations in 1D in the absence of convection on~$\IR$, under the condition that $\varphi(x)=\pm1$ as $x\to\pm\infty$:
\begin{equation}
\label{eq:CH1D}
\nabla\cdot{}m\nabla\big(\sigma\varepsilon\Delta\varphi^{\varepsilon}-\sigma\varepsilon^{-1}\Psi'(\varphi^{\varepsilon})\big)=0
\quad\text{on }\IR,
\quad\varphi^{\varepsilon}(x)=\pm{}1\quad\text{as }x\to\pm\infty \,,
\end{equation}
with the double-well potential~$\Psi$ according to~\eqref{eq:closureeqs_sub4}. A solution to~\eqref{eq:CH1D} is provided by $\varphi^{\varepsilon}(x)=\tanh(x/\sqrt{2}\varepsilon)$; see, e.g., \cite{Demont:2022dk} and also~\EQ{tanh} below. The solution can be conceived of as the phase field associated with a diffuse interface centered at $x=0$. Our analysis is based on the assumption that the phase-field generally exhibits the $\tanh(\cdot)$\nobreakdash-form across a fluid--fluid interface with respect to the normal to that interface, and that derivatives tangential to the interface are negligible. The latter implies that the curvature of the interface is negligible, in particular, $\kappa\varepsilon\ll{}1$ with $\kappa$ as the additive curvature of the interface. 

We denote by $\Gamma_{\LA}$ the zero level-set of~$\varphi$. We regard the Ginzburg--Landau energy density in~\eqref{eq:nrg_functnl}. The Ginzburg--Landau-energy density is non-negative and, under the above assumptions, its mass concentrates on~$\Gamma_{\LA}$ as $\varepsilon\to{}+0$. Denoting by $\psi$ a suitably smooth function, it  holds that
\begin{multline}
\label{eq:GL2ST}
\int_{\Omega}\psi\,\left( \frac{\sigma\varepsilon}{2} | \nabla \varphi |^2 + \frac{\sigma}{\varepsilon} \Psi \right)
=
\int_{\Gamma_{\LA}}
\int_{-a_{\varepsilon}}^{a_{\varepsilon}}
\psi(\B{x}+s\B{n}_{\LA})
\left( \frac{\sigma\varepsilon}{2} (d_s\varphi^{\varepsilon}(s))^2 + \frac{\sigma}{\varepsilon} \Psi(\varphi^{\varepsilon}(s)) \right)\,\dd{}s\,\dd\Gamma(\B{x})
\\
=
\int_{\Gamma_{\LA}}
\int_{-a_{\varepsilon}}^{a_{\varepsilon}}
\psi(\B{x}+s\B{n}_{\LA})\,\frac{2\sqrt{2}\sigma}{3}\delta_{\varepsilon}(s)\,\dd{}s\,\dd\Gamma(\B{x})
=
\sigma_{\LA}
\int_{\Gamma_{\LA}}
\psi(\B{x})\,\dd\Gamma(\B{x}) \,,
\end{multline}
as $\varepsilon\to{}+0$, for any positive $\varepsilon\ll{}a_{\varepsilon}\ll{}1$ and with
\begin{equation}
\delta_{\varepsilon}(s)=\frac{3}{4\sqrt{2}\varepsilon}\operatorname{sech}^4\bigg(\frac{s}{\sqrt{2}\varepsilon}\bigg)     \,.
\end{equation}
The vector $\B{n}_{\LA}$ denotes the unit normal vector on~$\Gamma_{\LA}$ exterior to the ambient domain. The function $\delta_{\varepsilon}(\cdot)$ can be conceived of as a regularized delta distribution. Equation~\eqref{eq:GL2ST} conveys that the Ginzburg--Landau energy density converges to a delta distribution on~$\Gamma_{\LA}$, multiplied by liquid-ambient surface tension, as $\varepsilon\to{}+0$.

\begin{figure}[!t]
\centering
\includegraphics[width=0.9\textwidth]{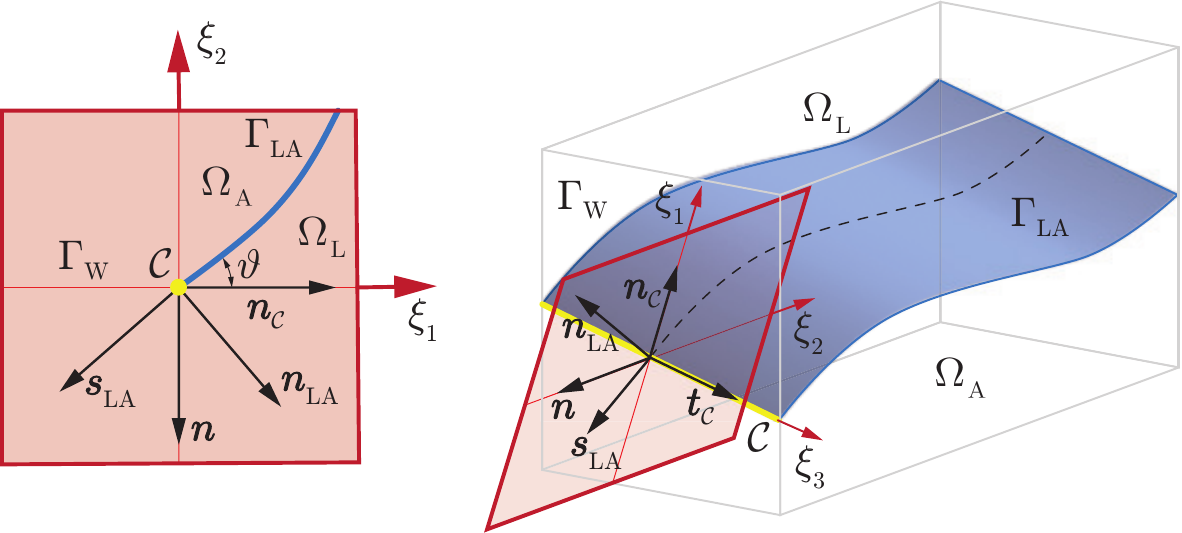}
\caption{Schematic of the considered contact-line configuration pertaining to the sharp-interface GNBC~(\ref{eq:SI-GNBC}), including nomenclature. The $(\xi_1,\xi_2)$ plane is orthogonal to the contact line. Adapted from~\cite{Gerbeau:2009km}.
\label{fig:GNBC}}
\end{figure}

We next consider the limit of the GNBC; see the illustration in Figure~\FIG{GNBC} for nomenclature. The interface $\Gamma_{\LA}$ and wetting boundary $\Gamma_{\W}$ intersect at the contact line~$\mathcal{C}$ at angle~$\vartheta$, interior to the liquid domain. 
The sharp-interface GNBC can be expressed as~\cite{Gerbeau:2009km}:
\begin{equation}
\label{eq:SI-GNBC}
\B{u}_S-\B{u}_{\W}=
\nu_2\big((p\B{n}-\B{\tau}\B{n})_S
-\sigma_{\LA}\big(\cos(\vartheta_{\textsc{eq}})+\B{s}_{\LA}\cdot\B{n}_{\mathcal{C}}
\big)\B{n}_{\mathcal{C}}\delta_{\mathcal{C}}\big) \,,
\end{equation}
where $\delta_{\mathcal{C}}$ is a delta distribution on the contact line, $\B{s}_{\LA}$ is the co-normal to~$\partial\Gamma_{\LA}$ and~$\B{n}_{\mathcal{C}}$ is the in-plane (parallel to~$\Gamma_{\W}$) unit normal vector on the contact line, oriented toward the liquid domain, and~$\vartheta_{\textsc{eq}}$ is the equilibrium contact angle. We first examine the correspondence between~\eqref{eq:wetting_BC_sub2} and~\EQ{SI-GNBC} tangential to the contact line. By virtue of the orthogonality of $\B{n}_{\mathcal{C}}$ and~$\B{n}$ to~$\B{t}_{\mathcal{C}}$, the sharp-interface GNBC~\EQ{SI-GNBC} implies
\begin{equation}
 \B{t}_{\mathcal{C}}\cdot(\B{u}_S-\B{u}_{\W})+\B{t}_{\mathcal{C}}\cdot\nu_2\B{\tau}\B{n}=0 \,.
\end{equation}
The diffuse-interface GNBC~\eqref{eq:wetting_BC_sub2} yields
\begin{equation}
\label{eq:GNBCt}
 \B{t}_{\mathcal{C}}\cdot(\B{u}_S-\B{u}_{\W})+\B{t}_{\mathcal{C}}\cdot\nu_2\B{\tau}\B{n}=
 \B{t}_{\mathcal{C}}\cdot\sigma_{\SF}'(\varphi)\nabla_S\varphi-\B{t}_{\mathcal{C}}\cdot\nu_2\B{\zeta}\B{n} \,.
\end{equation}
Because $\nabla\varphi$ is oriented along $\B{n}_{\LA}$ by virtue of the standing assumption on the phase field, it is orthogonal to~$\B{t}_{\mathcal{C}}$, and the first right-hand-side term in~\eqref{eq:GNBCt} vanishes. For the second term, it follows from~\eqref{eq:closureeqs_sub3} that
\begin{equation}
 \B{t}_{\mathcal{C}}\cdot\nu_2\B{\zeta}\B{n}
 =
 \nu_2\Big(\B{t}_{\mathcal{C}}\cdot\B{n}\left( \frac{ \sigma \varepsilon } { 2 } | \nabla \varphi |^2 + \frac{ \sigma } { \varepsilon } \Psi \right)
 - \B{t}_{\mathcal{C}}\cdot\sigma \varepsilon \nabla_S \varphi\, \partial_n\varphi\Big)   \,,
\end{equation}
which again vanishes on account of orthogonality. Hence, the diffuse-interface GNBC and sharp-interface GNBC are consistent in the direction tangential to the contact line. To compare the diffuse-interface GNBC and sharp-interface GNBC in the direction~$\B{n}_{\mathcal{C}}$, we note that the phase field in the vicinity of the contact line can be expressed as
\begin{equation}
\label{eq:varphiC}
\hat{\varphi}(\xi_1,\xi_2)=\tanh\bigg(\frac{\xi_1\sin(\vartheta)-\xi_2\cos(\vartheta)}{\sqrt{2}\varepsilon}\bigg)     \,,
\end{equation}
with reference to the $\B{\xi}:=(\xi_1,\xi_2,\xi_3)$ coordinates in Figure~\ref{fig:GNBC}. Assuming that the contact line moves in the direction~$\B{n}_{\mathcal{C}}$ with velocity~$v$, the local coordinates~($\B{\xi}$) and global~($\B{x}$) coordinates are related by $\B{x}=(\xi_1+vt)\B{n}_{\mathcal{C}}-\xi_2\B{n}+\xi_3\B{t}_{\mathcal{C}}$. For the $\varphi$\nobreakdash-dependent terms in the right-hand side of~\eqref{eq:wetting_BC_sub2}, we then obtain
\begin{equation}
\begin{aligned}
&\big[\B{n}_{\mathcal{C}}\cdot\nu_2\big(-\B{\zeta}\B{n}+  \nabla_S \sigma_{\SF}(\varphi)\big)\big](t,\B{x})
\\
&\qquad=
\nu_2\big[\sigma\varepsilon\partial_{\xi_1}
\hat{\varphi}\,\partial_{\xi_2}\hat{\varphi}+\sigma_{\SF}'(\hat{\varphi})\partial_{\xi_1}\hat{\varphi}\big]\big(\B{x}\cdot\B{n}_{\mathcal{C}}-vt,0,\B{x}\cdot\B{t}_{\mathcal{C}}\big)
\\
&\qquad=\big[\sigma_{\SL}-\sigma_{\SA}+\sigma_{\LA}\cos(\vartheta)\big]\,\sin(\vartheta)\,\delta_{\varepsilon}\big((\B{x}\cdot\B{n}_{\mathcal{C}}-vt)\sin(\vartheta)\big)
\\
&\qquad=
-\sigma_{\LA}[\cos(\vartheta_{\textsc{eq}})-\cos(\vartheta)]\,\delta_{\varepsilon}^{\star}(\B{x}\cdot\B{n}_{\mathcal{C}}-vt) \,,
\end{aligned}
\end{equation}
where $\delta_{\varepsilon}^{\star}(\cdot)=\sin(\vartheta)\delta_{\varepsilon}((\cdot)\sin(\vartheta))$ represents a regularized delta distribution. In the ultimate identity, we used the relation between the equilibrium contact angle, $\vartheta_{\textsc{eq}}$, and the liquid-ambient, solid-ambient, and solid-liquid surface tensions, viz. $\sigma_{\LA}\cos\vartheta_{\textsc{eq}}=\sigma_{\SA}-\sigma_{\SL}$. Noting that $\B{s}_{\LA}\cdot\B{n}_{\mathcal{C}}=-\cos(\vartheta)$, it follows that the diffuse-interface GNBC~\eqref{eq:wetting_BC_sub2} indeed formally converges to the sharp-interface GNBC~\eqref{eq:SI-GNBC} as~$\varepsilon\to+0$.

\subsection{Contact-angle condition}
\label{sec:CA_condition}
The wetting-boundary condition~\eqref{eq:wetting_BC_sub1} can be conceived of as a dynamic contact-angle condition. In the present work, we concentrate on the effect of slip, which is encoded in the GNBC, and we focus on the equilibrium contact-angle condition that emerges in the limit~$\nu_1\to\infty$, viz. the nonlinear Robin-type condition
\begin{equation}
\label{eq:eqconang}
\sigma\varepsilon\partial_n\varphi+\sigma_{\SF}'(\varphi)=0,
\qquad\text{on }\Gamma_{\W}\,.
\end{equation}
In this section, we establish that if $\vartheta=\vartheta_{\textsc{eq}}$, \cref{eq:varphiC} is compatible with~\eqref{eq:eqconang}, and with~\eqref{eq:wetting_BC_sub1} if the contact-line velocity~$v$  coincides with the flow velocity orthogonal to the contact line, $\B{u}_S\cdot\B{n}_{\mathcal{C}}$. The former has also been derived in~\cite{Shokrpour-Roudbari:2016dp}, and we repeat it here for completeness and as an intermediate step towards the latter. 

Proceeding formally, we substitute~\eqref{eq:varphiC} in~\eqref{eq:eqconang}:
\begin{equation}
\label{eq:eqconang1}
\begin{aligned}
\big[\sigma\varepsilon\partial_n\varphi+\sigma_{\SF}'(\varphi)\big](t,\B{x})
&=
\Big[
-\sigma\varepsilon\partial_{\xi_2}\hat{\varphi}
+\frac{3}{4}(\hat{\varphi}{}^2-1)(\sigma_{\SA}-\sigma_{\SL})\Big]\big(\B{x}\cdot\B{n}_{\mathcal{C}}-vt,0,\B{x}\cdot\B{t}_{\mathcal{C}}\big)
\\
&=
\frac{3}{4}\sigma_{\LA}\big(\cos(\vartheta)-\cos(\vartheta_{\textsc{eq}})\big)
\operatorname{sech}^2\bigg(\frac{(\B{x}\cdot\B{n}_{\mathcal{C}}-vt)\sin(\vartheta)}{\sqrt{2}\varepsilon}\bigg) \,.
\end{aligned}
\end{equation}
The right-hand side of~\eqref{eq:eqconang1} indeed vanishes for $\vartheta=\vartheta_{\textsc{eq}}$, independent of~$\varepsilon$.

Considering the left-hand side of the dynamic contact-angle condition~\eqref{eq:wetting_BC_sub1}, substitution of~\eqref{eq:varphiC} yields
\begin{equation}
\big[\partial_t\varphi
+
\B{u}_{S}\cdot\nabla_{S}\varphi\big](t,\B{x})
=
(\B{u}_{S}\cdot\B{n}_{\mathcal{C}}-v)\,
\frac{\sin(\vartheta)}{\sqrt{2}\varepsilon}\operatorname{sech}^2\bigg(\frac{(\B{x}\cdot\B{n}_{\mathcal{C}}-vt)\sin(\vartheta)}{\sqrt{2}\varepsilon}\bigg)\,,
\end{equation}
which indeed vanishes if~$v=\B{u}_{S}\cdot\B{n}_{\mathcal{C}}$, independent of $\varepsilon$. In combination with~\eqref{eq:eqconang1}, this implies that the diffuse-interface profile~\eqref{eq:varphiC} is compatible with~\eqref{eq:wetting_BC_sub1} if the contact line moves with the flow velocity.


\section{Numerical experiments}
\label{sec:num_exp}

In this section, we illustrate several properties of the diffuse-interface NSCH model~\EQ{strong} by means of a Couette flow problem. In~\cref{ssec:numexp_noslip}, we show that the diffuse-interface model does not converge as $m \propto \varepsilon^2 \rightarrow +0$ for the case of no-slip wetting boundary conditions, which is in line with the velocity-field inconsistency at the contact point in the classical sharp-interface free-boundary problem as presented in \cite{huh1971hydrodynamic}. Furthermore, we measure, for several values of both the interface thickness $\varepsilon$ and the \Ca[l] \Ca, the deviation of a fixed-mobility no-slip diffuse-interface model to the classical free-boundary problem, in which the inherent mobility-induced slip length of the diffuse-interface model corresponds to the GNBC slip length of the sharp-interface model~\cite{Jacqmin:2000kx}. Next, in~\cref{ssec:numexp_GNBC}, we introduce the wetting and generalized Navier boundary conditions to the Abels--Garcke--Gr\"{u}n NSCH model. We investigate, when introducing finite slip, the convergence to the sharp-interface free-boundary problem as $\varepsilon, m \rightarrow +0$ for various values of the slip length. Finally, in~\cref{sssec:2BC}, we investigate the formation of the so-called `triple-wedge flows' in the case of non-uniform viscosity ratios, and separately for the case of a non-right-angle microscopic contact angles.

For our diffuse-interface simulations, we make use of the extensively verified higher-order finite element code described in \cite{Demont_Stoter_van_Brummelen_2023,Demont:2022dk,Brummelen:2021aw}\footnote{All diffuse-interface simulations have been performed using the open source software package \mbox{Nutils~(\url{https://www.nutils.org}).}}. The sharp-interface results are obtained with the finite element code described in \cite{diddens2023bifurcation}\footnote{All sharp-interface simulations have been performed using the open source software package \mbox{Pyoomph~(\url{https://pyoomph.github.io}).}}. Where necessary, we employ (up to eight levels of) spatial and temporal adaptive refinement to ensure convergence of all presented results.


\subsection{Set-up}
\label{sec:setup}
\begin{figure}[!t]
\centering
\includegraphics[width=0.90\textwidth]{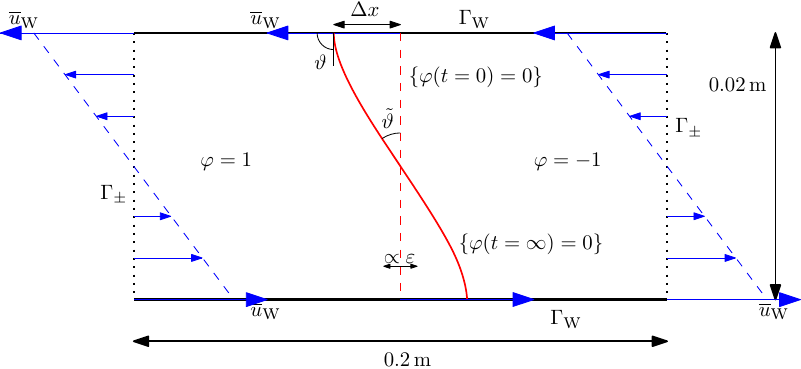}
\caption{Schematic overview of the considered Couette test case. The dashed red line indicates the $\{\varphi(t=0) = 0\}$ contour line of the initial phase-field. The solid red line illustrates the equilibrated $\{\varphi(t=\infty)=0\}$ contour, representing the interface between the two fluids, corresponding to wall 
velocity~$\overline{u}_\W$ at the wetting boundaries~$\Gamma_{\W}$. Quantities of interest are the contact-point displacement, $\Delta x$, and macroscopic contact angle,~$\tilde{\vartheta}$.}
\label{fig:Couette_model_IC_BC}
\end{figure}

The physical domain of the considered Couette flow system is given by $\Omega = (0, 0.2)\,\text{m} \times (0, 0.02)\,\text{m}$. The chosen $10:1$ domain length ratio minimizes the effects of the in- and outflow boundary conditions. A schematic overview of the set-up for the test cases is presented in Figure \FIG{Couette_model_IC_BC}.

On the top and bottom boundaries, $\Gamma_{\W} \coloneqq (0, 0.2)\,\text{m} \times \{0, 0.02\}\,\text{m}$, we prescribe the static contact angle boundary condition~\eqref{eq:eqconang} and the generalized Navier boundary conditions~\eqref{eq:wetting_BC_sub2}, complemented by the impermeability condition~\eqref{eq:wetting_BC_sub3} and the homogenous Neumann condition for the chemical potential~\eqref{eq:wetting_BC_sub4}. On the left and right boundaries, $\Gamma_{\pm} \coloneqq \partial \Omega \setminus \Gamma_{\W}$, we prescribe pre-defined Dirichlet in- and outflow conditions for the velocity field, and homogeneous Neumann conditions for the order parameter and the chemical potential, according to~\EQ{pm_BC}. Here, the $x_2$-component of the function $\overline{\B{u}}_{\pm}(x_2,t)$ is set to zero, and the $x_1$-component is taken as a linear profile with odd symmetry about the centerline. Its maximum velocity at the top and bottom walls is proportional to the traction and slip, i.e.
\begin{equation}
\label{eq:wall_vel_eq}
(\overline{u}_1)_{\pm}\big|_{x_2=0} - u_{\W} = \nu_2 \, \eta \, \partial_{x_2}(\overline{u}_1)_{\pm}\big|_{x_2=0} \, ,
\end{equation}
where $u_{\W}$ is a given external wall velocity according to~\eqref{eq:uW} below. This boundary condition conforms with the single-species laminar Couette flow profile. The time-dependent wall velocity $u_\W$ and, thus, also the in-/outflow velocity $(\overline{u}_1)_{\pm}$ in~\eqref{eq:wall_vel_eq}, are ramped up over the course 
of~$1$ second according to a sinusoidal profile with zero derivative at its minimum and maximum values:
\begin{equation}
\label{eq:uW}
u_\W(t) = \left\{
\begin{tabular}{ll}
$(\frac{1}{2} - \frac{1}{2}\cos(\pi t)) \overline{u}_\W$, & $t<1\,\text{s}$\,,\\
$\overline{u}_\W$, & $t \geq 1\,\text{s}$\,,
\end{tabular}
\right.
\end{equation}
with $\overline{u}_\W$ the maximum wall velocity.

We initialize the simulations by means of a stationary velocity field,
\begin{equation}
\B{u}(0,\B{x}) = \B{u}_0(\B{x}) \equiv 0\,.
\end{equation}
The initial order parameter $\varphi$ is such that it corresponds to a configuration in which the two fluids $\LL$ and $\AA$ are separated by a centered vertical fluid--fluid interface, $\Gamma_{\LA}$, making a $\pi/2$ microscopic contact angle with the wetting boundaries $\Gamma_{\W}$:
\begin{equation}
\label{eq:tanh}
\varphi(0,\B{x}) = \varphi_0(\B{x}) \coloneqq \tanh \left(\frac{d_{\pm}(\B{x},\Gamma_{\LA})}{\sqrt{2}\varepsilon}\right) = \tanh \left( \frac{-x_1 + 0.1 \,\text{m}}{\sqrt{2}\varepsilon} \right),
\end{equation}
where $d_{\pm}(\B{x}, \Gamma_{\LA})$ represents the signed distance function from $\B{x}$ to $\Gamma_{\LA}$. The initial phase field~$\varphi_0(\B{x})$ coincides with an equilibrium solution of the Cahn--Hilliard equations for the phase field in one spatial dimension (see \cref{sec:SI_lim_analysis}), and attributes equal volumes to both fluids in~$\Omega$.

The simulations are considered equilibrated when further changes in the solution fields no longer have a noticeable effect on the results. The default parameter values used in all numerical experiments throughout this section are reported in 
Table~\TAB{model_param_vals}.

\begin{table}[!t]\centering
\caption{Physical and model parameter values for the diffuse-interface NSCH model of the considered numerical experiments. Entries marked with the symbol~$\ast$ or in-between accolades  indicate a set of values, which will be (further) specified in the text. Entries in-between brackets indicate properties of the individual $\LL$ and $\AA$ phases of the fluid mixture.}
\label{tab:model_param_vals}
\resizebox{\textwidth}{!}{
\begin{tabular}{c|cc|ccc|cccc|cc}
&\multicolumn{2}{c|}{Fluids}&\multicolumn{3}{c|}{Interface}&\multicolumn{4}{c|}{Boundary}&\multicolumn{2}{c}{Slip Lengths}\\
TC\# & $\rho_{\LL}$, $\rho_{\AA}$ & $\eta_{\LL}$, $\eta_{\AA}$ & $\sigma_{\LA}$  & $\varepsilon$ & $m$ & $\nu_1$  & $\nu_2$ & $\overline{u}_\W$ & $\vartheta$ & $s_m$ & $s_\nu$ \\[0.1cm]
\hline
\hline & & & & & & & & & & & \\[-0.3cm]
--- & $\frac{\textrm{kg}}{\textrm{m}^d}$ &$\frac{\textrm{kg}\,\textrm{m}^{2-d}}{\textrm{s}}$& $\frac{\textrm{kg}\,\textrm{m}^{3-d}}{\textrm{s}^2}$ & \textmu$\text{m}$ & $\frac{\textrm{m}^d\,\textrm{s}}{\textrm{kg}}$ & $\frac{\text{m}^{d-3}\,\text{s}}{\text{kg}}$ & $\frac{\textrm{m}^{d-1}\,\textrm{s}}{\textrm{kg}}$ & $\frac{\textrm{m}}{\textrm{s}}$ & $\text{rad}$ &  $\text{mm}$ & $\text{mm}$ \\[0.2cm]
1A & $10^3$ & $10^{-1}$ & $7.28\!\times\!10^{-2}$ & $\ast$ & $\{4,1\}\!\times\!10^{-5}$ & $\infty$ & $0$ & $4\!\times\!10^{-3}$ & $\pi/2$ & $\{2,1\}$ & $0$ \\[0.2cm]
1B & $10^3$ & $10^{-1}$ & $7.28\!\times\!10^{-2}$ & $50$ & $\{4,1\}\!\times\!10^{-5}$ & $\infty$ & $0$ & $\ast$ & $\pi/2$ & $\{2,1\}$ & $0$ \\[0.2cm]
1C & $10^3$ & $10^{-1}$ & $7.28\!\times\!10^{-2}$ & $\ast$ & $\propto \varepsilon^2$ & $\infty$ & $0$ & $4\!\times\!10^{-3}$ & $\pi/2$ & $\propto \varepsilon$ & $0$ \\[0.2cm]
\hline & & & & & & & & & & & \\[-0.3cm]
2A & $10^3$ & $10^{-1}$ & $7.28\!\times\!10^{-2}$ & $\ast$ & $\propto \varepsilon^2$ & $\infty$ & $\{2,1\}\!\times\!10^{-2}$ & $4\!\times\!10^{-3}$ & $\pi/2$ & $\propto \varepsilon$ & $\{2,1\}$ \\[0.2cm]
2B & $10^3$ & $10^{(-1,-3)}$ & $7.28\!\times\!10^{-2}$ & $50$ & $10^{-9}$ & $\infty$ & $1\!\times\!10^{-2}$ & $4\!\times\!10^{-3}$ & $\pi/2$ & $10^{(-2,-3)}$ & $(1, 10^{-2})$ \\[0.2cm]
2C & $10^3$ & $10^{-1}$ & $7.28\!\times\!10^{-2}$ & $50$ & $10^{-9}$ & $\infty$ & $1\!\times\!10^{-2}$ & $4\!\times\!10^{-3}$ & $\pi/4$ & $10^{-2}$ & $1$ \\[0.1cm]
\end{tabular}
}
\end{table}


\subsection{Classical no-slip results}
\label{ssec:numexp_noslip}

\begin{table}[!b]\centering
\caption{Test case 1A: contact point displacements $\Delta x$, macroscopic contact angle $\tilde{\vartheta}$, and shear force~$F_\textsc{S}$ for the diffuse-interface model with no-slip boundary conditions, for slip lengths $s_{\nu,m} = \{2,1\}\,\textrm{mm}$, a decreasing sequence of diffuse-interface-thickness parameters~$\varepsilon$, and fixed mobility $m=\{4,1\}\!\times\!10^{-5}\,\textrm{m}^d\,\textrm{s}\,\textrm{kg}^{-1}$. Bottom row displays sharp-interface results with slip boundary conditions.}
\label{tab:1A_m_fixed_comb}
\begin{tabular}{c|ccc|ccc}
 & \multicolumn{3}{c|}{$s_{\nu,m}=2\,\text{mm}$} & \multicolumn{3}{c}{$s_{\nu,m}=1\,\text{mm}$} \\[0.2cm]
\hline & & & & & & \\[-0.3cm]
$\varepsilon$ & $\Delta x$ & $\tilde{\vartheta}$ & $F_\textsc{S}$ & $\Delta x$ & $\tilde{\vartheta}$ & $F_\textsc{S}$ \\[0.1cm]
\hline
\hline & & & & & & \\[-0.3cm]
$1.6\!\times\!10^{-3}\,\text{m}$ & $10^{-4}\,\text{m}$ & $10^{-2}\,\text{rad}$ & $10^{-3}\,\frac{\text{kg}\,\text{m}}{\text{s}^2}$ & $10^{-4}\,\text{m}$ & $10^{-2}\,\text{rad}$ & $10^{-3}\,\frac{\text{kg}\,\text{m}}{\text{s}^2}$ \\[0.2cm]
$2^{0}\phantom{^{-}}$ & $5.101$ & $7.602$ & $3.078$ & $9.098$ & $12.56$ & $5.469$ \\[0.2cm]
$2^{-1}$              & $5.380$ & $7.881$ & $3.228$ & $9.226$ & $12.43$ & $5.516$ \\[0.2cm]
$2^{-2}$              & $5.566$ & $8.060$ & $3.334$ & $9.407$ & $12.52$ & $5.613$ \\[0.2cm]
$2^{-3}$              & $5.685$ & $8.176$ & $3.403$ & $9.568$ & $12.65$ & $5.704$ \\[0.2cm]
$2^{-4}$              & $5.756$ & $8.248$ & $3.445$ & $9.680$ & $12.75$ & $5.770$ \\[0.2cm]
$2^{-5}$              & $5.795$ & $8.288$ & $3.468$ & $9.748$ & $12.81$ & $5.810$ \\[0.2cm]
\hline & & & & & & \\[-0.3cm]
NS w/slip             & $6.171$ & $7.836$ & $3.077$ & $8.848$ & $10.94$ & $4.801$ \\[0.1cm]
\end{tabular}
\end{table}
To illustrate the properties of the no-slip NSCH model, we perform several sets of numerical experiments. For these experiments, we enforce the fluid--fluid interface to make a $\pi/2 \, \text{rad}$ microscopic contact angle with the solid boundary, and impose no-slip boundary conditions on the wetting boundaries $\Gamma_\W$ by setting $\nu_2 = 0 \, \textrm{kg}^{-1} \, \textrm{m}^{d-1} \, \textrm{s}$. As such, any slip present in the diffuse-interface model is induced by the interfacial equilibration controlled by the mobility parameter $m$ as shown in~\cite{Jacqmin:2000kx}. We measure several quantities of interest. The contact point displacements $\Delta x$, measured relative to the initially imposed contact points at $x = 0.1\,\text{m}$, are determined from the phase-field $\varphi = 0$ contact points on $\Gamma_\W$. The macroscopic contact angle deviation $\tilde{\vartheta}$, or mid-box angle deviation, is determined from the direction of~$\nabla \varphi$ at the mid-box $\B{x} = (0.1, 0.01)\,\text{m}$ relative to the direction vector $(-1,0)$. Finally, we compute the excess shear force $F_\textsc{S}$ on the wetting boundaries $\Gamma_{\W}$ according to
\begin{equation}
F_\textsc{S} = - \int_{\Gamma_{\W}} \eta \partial_{x_2} u_{x_1}\,\text{d}x + \sigma \varepsilon \int_{\Gamma_{\W}} \partial_{x_1}\varphi \partial_{x_2}\varphi\,\text{d}x + \int_{\Gamma_{\W}} \eta \, \partial_{x_2}(\overline{u}_1)_{\pm}\big|_{x_2=0}\, \text{d}x\,.
\end{equation}
where the second term vanishes for $\vartheta=\pi/2$, and the latter term cancels the inherent shear force acting on the wetting boundaries due to finite distance between both wetting boundaries in the Couette flow problems. We subtract this `baseline' contribution in the shear force computation in order to emphasize the contact point contributions.

We compare the results of each set of numerical experiments with the results of simulations with the classical sharp-interface counterpart, on which slip boundary conditions are imposed with the corresponding slip length parameter \cite{Jacqmin:2000kx}: we equate the GNBC slip length, $s_\nu \coloneqq \eta \nu_2$, of the sharp-interface simulations to the mobility-induced slip length, $s_m \coloneqq \sqrt{\eta m}$, of the diffuse-interface simulations, by choosing appropriate values of $\nu_2$ and $m$.

\subsubsection{Test case 1A: NS/MS limit}
\label{sec:1A}

We first consider the case where we let $\varepsilon$ approach zero, while keeping the mobility $m$, and consequently the mobility-induced slip length, $s_m$, fixed. Two separate slip lengths, $s_m=2\,\text{mm}$ and $s_m=1\,\text{mm}$, are considered; see test case 1A in Table \TAB{model_param_vals}. We compare the diffuse-interface results with the results from sharp-interface simulations in which the sharp-interface slip lengths on the wetting boundaries $\Gamma_\W$, $s_\nu$, are matched with the mobility induced diffuse-interface slip lengths, $s_m$. As $m>0$ in this test case is fixed, the NSCH model converges to the non-classical Navier--Stokes/Mullins--Sekerka model as $\varepsilon \rightarrow +0$, as established in~\cite[\S{}4.1]{Abels:2018ly}. Thus, a deviation in results from the classical sharp-interface model with slip is expected when approaching the limit.

\begin{figure}[!b]
\centering
\includegraphics[width=1.00\textwidth]{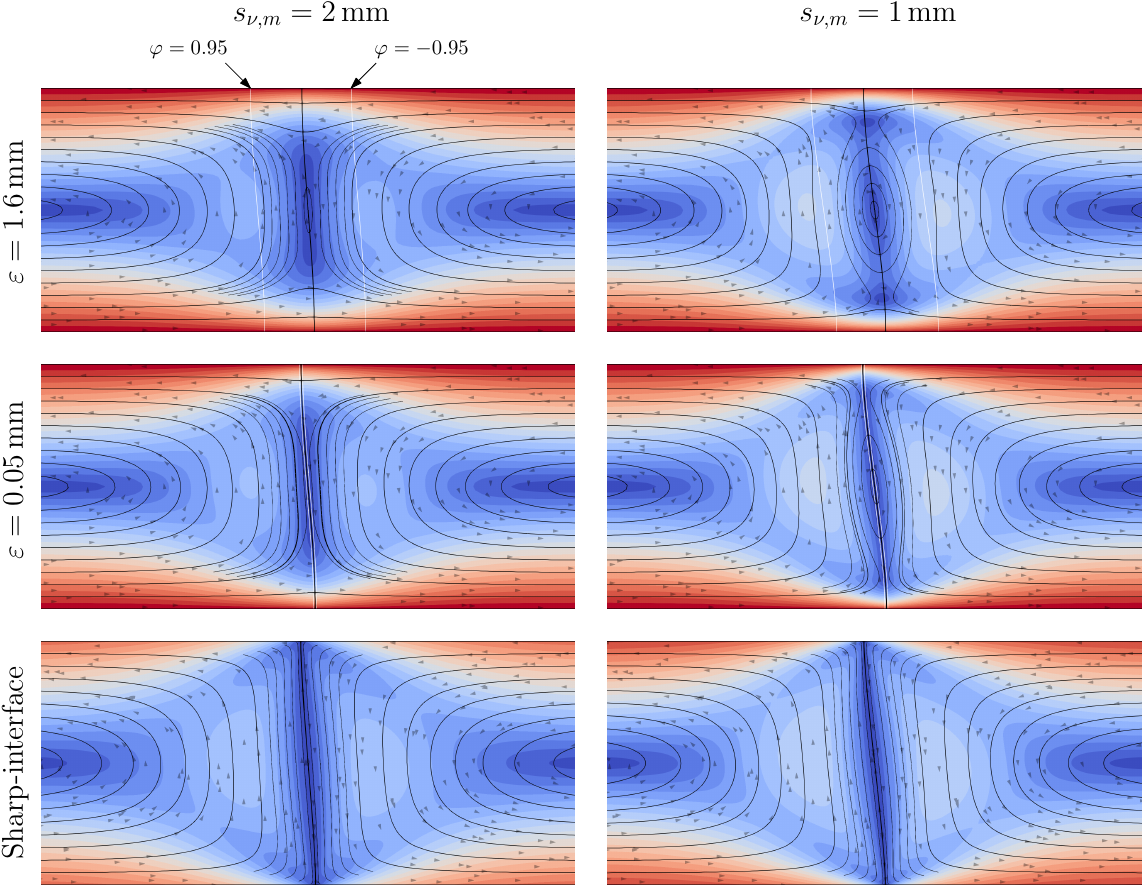}
\caption{No-slip NSCH (top) model with different mobility-induced slip lengths ($s_m$) and interface widths ($\varepsilon$), and sharp-interface model with different slip lengths ($s_\nu$). Showing velocity fields with streamlines.}
\label{fig:1A_velocity_fields}
\end{figure}

The quantities of interest are displayed in Table \TAB{1A_m_fixed_comb}. From the successive values of $\Delta x$, $\tilde{\vartheta}$, and $F_\textsc{S}$ for decreasing values of $\varepsilon$, it is clear that for both slip lengths, the diffuse interface model does not converge to the classical sharp-interface model. While the results of the sharp interface model with non-zero slip deviate slightly from the limiting results of the diffuse interface simulations as $\varepsilon \rightarrow +0$, it is noteworthy that the results do roughly coincide for both slip lengths, substantiating the asymptotic results in \cite{Jacqmin:2000kx}. However, the limiting contact point displacement of the diffuse-interface model increases much faster than the contact point displacement of the sharp-interface model as the slip length decreases, suggesting a larger critical slip length for the NSCH model where the contact points drift off to $\pm \infty$. Likewise, the limiting shear force of the diffuse-interface model grows faster than for the sharp-interface model. 

Interestingly, for the larger slip length in Table \TAB{1A_m_fixed_comb}, the contact point displacement $\Delta x$ of the NSCH model is smaller than that of the sharp-interface NS model, while the macroscopic contact angle deviation $\tilde{\vartheta}$ is larger, implying the diffuse-interface model exhibits larger interface curvature and larger transitional zones between the microscopic and macroscopic contact angle regions. Likewise, for the smaller slip length, the ratio $\tilde{\vartheta} / \Delta x$ is larger for the NSCH model than for the NS-with-slip model, resulting in the same qualitative conclusions about the interface shape. This mismatch is indicative of the fact that the constant-mobility sharp-interface limit of the NSCH model does not yield the classical binary-flow Navier-Stokes model, but the Mullins--Sekerka model instead.

In Figure~\FIG{1A_velocity_fields}, the velocity fields of the diffuse-interface and sharp-interface simulations are presented for several slip lengths and values of $\varepsilon$. In the diffuse-interface results one observes a positive flux through the fluid--fluid interface near the wetting boundaries $\Gamma_\W$, which remains positive as $\varepsilon \rightarrow 0$. Prominently, the fluid circulation near the center of --- and through --- the interface, grows stronger as the slip length $s_m$ becomes smaller. Additionally, the presence of this circulation through the interface for $\varepsilon = 0.05\,\textrm{mm}$ and $s_m = 1\,\textrm{mm}$, where the dimensions of the recirculation zone are much larger than the interface thickness, indicates that it is part of the limiting solution and not invoked by the finite interface thickness. Naturally, in the sharp-interface counterpart, there is no flux through the fluid--fluid interface, and the fluid velocity magnitude is smaller along $\Gamma_\W$ in its entirety, vanishing at the contact points. As such, the flow profiles of both models are significantly different. The disparity in velocity magnitudes away from the contact point is ${O}(s_\cdot)$ as $s_\cdot \rightarrow 0$. It is to be noted, however, that the no-slip limit itself does not exist, as finite displacement of the contact point in this limit would lead to a contradiction. These results convey that mobility-mediated slip in a diffuse-interface model, i.e. slip induced by non-vanishing mobility in combination with a no-slip condition on velocity, while providing approximately the same interface profile as the sharp-interface model, does not yield the same flow behavior in the vicinity of the interface.

\begin{table}
\centering
\caption{Test case 1B: contact point displacements $\Delta x$, macroscopic contact angle $\tilde{\vartheta}$, and shear force~$F_\textsc{S}$ for the diffuse-interface (DI) model with no-slip boundary conditions, for various values of the wall velocity $\overline{u}_\W$, for slip lengths $s_{\nu,m} = \{2,1\}\,\textrm{mm}$s and with 
fixed~$\varepsilon=5\times{}10^{-5}\,\text{m}$ 
and~$m=\{4,1\}\!\times\!10^{-5}\,\textrm{m}^d\,\textrm{s}\,\textrm{kg}^{-1}$. The sharp-interface (SI) results are with slip wetting boundary conditions.
\label{tab:uW_comb}}
\begin{tabular}{c|c|ccc|ccc}
& & \multicolumn{3}{c|}{NSCH (DI)} &\multicolumn{3}{c}{NS with slip (SI)} \\
$s_{\nu,m}$ & $\overline{u}_\W$ & $\Delta x$ & $\tilde{\vartheta}$ & $F_\textsc{S}$ & $\Delta x$ & $\tilde{\vartheta}$ & $F_\textsc{S}$ \\[0.1cm]
\hline
\hline & & & & & & & \\[-0.3cm]
$\text{mm}$ & $10^{-3}\,\frac{\text{m}}{\text{s}}$ & $10^{-4}\,\text{m}$ & $10^{-2}\,\text{rad}$ & $10^{-4}\,\frac{\text{kg}\,\text{m}}{\text{s}^{-2}}$ & $10^{-4}\,\text{m}$ & $10^{-2}\,\text{rad}$ & $10^{-4}\,\frac{\text{kg}\,\text{m}}{\text{s}^{-2}}$ \\[0.2cm]
\multirow{5}{*}{2} & $8$ & $11.75$ & $16.75$ & $69.76$ & $12.53$ & $15.86$ & $61.93$ \\[0.2cm]
& $4$ & $5.795$ & $8.288$ & $34.68$ & $6.171$ & $7.836$ & $30.77$ \\[0.2cm]
& $2$ & $2.887$ & $4.133$ & $17.31$ & $3.075$ & $3.907$ & $15.36$ \\[0.2cm]
& $1$ & $1.442$ & $2.057$ & $8.651$ & $1.536$ & $1.952$ & $7.678$ \\[0.2cm]
\hline & & & & & & & \\[-0.3cm]
\multirow{5}{*}{1} & $8$ & $20.26$ & $26.38$ & $118.1$ & $18.22$ & $22.39$ & $97.25$ \\[0.2cm]
& $4$ & $9.748$ & $12.81$ & $58.10$ & $8.848$ & $10.94$ & $48.01$ \\[0.2cm]
& $2$ & $4.831$ & $6.364$ & $28.94$ & $4.393$ & $5.439$ & $23.93$ \\[0.2cm]
& $1$ & $2.410$ & $3.177$ & $14.46$ & $2.193$ & $2.716$ & $11.96$ \\[0.1cm]
\end{tabular}
\end{table}

\subsubsection{Test case 1B: Capillary number}

\begin{figure}
\centering
\includegraphics[width=1.00\textwidth]{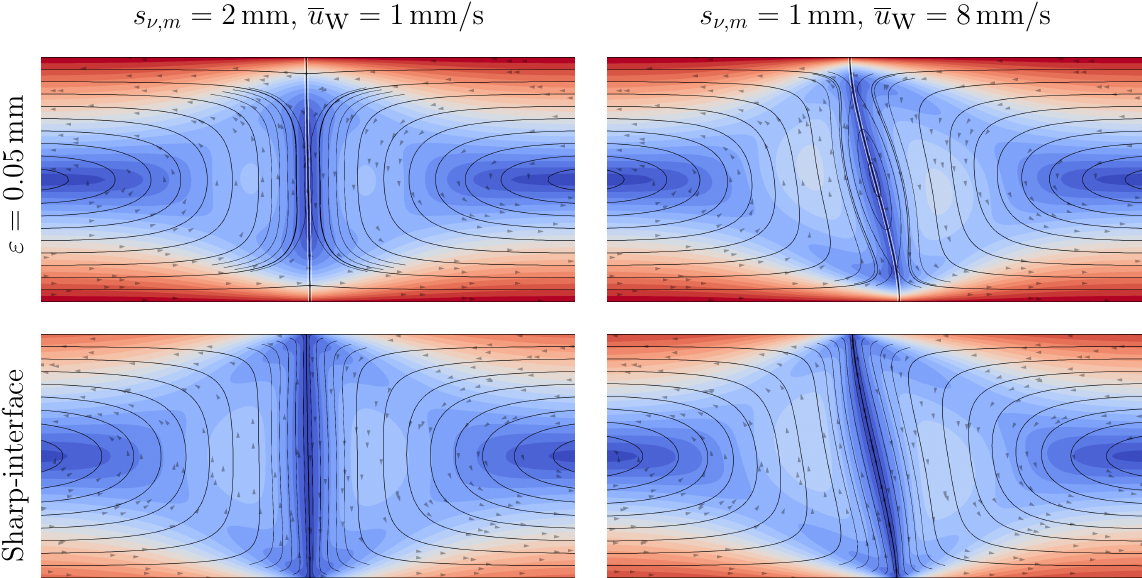}
\caption{No-slip NSCH (top) and sharp-interface (bottom) models with different slip lengths ($s_\nu,s_m$) and wall velocities ($\overline{u}_\W$). Showing velocity fields with streamlines.}
\label{fig:1B_velocity_fields}
\end{figure}

Next, let us investigate the dependence of the quantities of interest and flow profiles on the wall velocity $\overline{u}_\W$, and consequently, the \Ca[l] $\Ca \coloneqq \eta \overline{u}_\W / \sigma$. We consider NSCH simulations with fixed mobility values and a small interface thickness, $\varepsilon = 5\!\times\!10^{-5}\,\text{m}$, relative to the mobility-induced slip lengths, $s_m \in \{2,1\}\,\text{mm}$. Similarly, the sharp-interface simulations will be performed with the corresponding GNBC slip-length $s_\nu \in \{2,1\}\,\text{mm}$. The slip-lengths $s_m,s_\nu$ in turn are small relative to the macroscopic length scale provided by the channel semi-height, $\ell = 1\,\text{cm}$, to minimize effects of the finite channel height. The wall velocities under consideration are $u_\W \in \{8,4,2,1\}\,\text{mm}\,\text{s}^{-1}$; the remaining physical and model parameter values are listed under test case 1B in Table~\TAB{model_param_vals}. The computed equilibrated quantities are presented in Table~\TAB{uW_comb}, and Figure \FIG{1B_velocity_fields} displays the velocity fields of the performed simulations that have the smallest and largest contact point displacements.

Due to the small angular displacements considered in this test case, $\Delta x$, $\tilde{\vartheta}$, and $F_\textsc{S}$ behave essentially linearly in the \Ca[l] for both the diffuse- and sharp-interface models and for both slip lengths. Increasing the \Ca[l] further quickly results in non-linear behavior and, beyond a slip-length-dependent critical wall velocity, ultimately leads to the contact points escaping to $\pm \infty$. Halving the slip length causes an increase in all of the observed quantities, for both models and for all considered wall velocities. For each quantity of interest and both slip lengths, 
the ratio between the diffuse-interface result and the sharp-interface result is essentially constant, independent of $\Ca$. This indicates consistency of both models across a wide range of \Ca[l]s.

\begin{table}[!b]\centering
\caption{Test case 1C: 
contact point displacements $\Delta x$, macroscopic contact angle $\tilde{\vartheta}$, and shear force~$F_\textsc{S}$ for the diffuse-interface model with no-slip boundary conditions, for a decreasing sequence of interface-thickness parameters~$\varepsilon$ and mobility~$m \propto \varepsilon^2$.}
\label{tab:1C}
\begin{tabular}{ccc|ccc}
$\varepsilon$ & $m$ & $s_m$ & $\Delta x$ & $\tilde{\vartheta}$ & $F_\textsc{S}$ \\[0.1cm]
\hline
\hline & & & & & \\[-0.3cm]
$1.6\!\times\!10^{-3}\,\text{m}$ & $10^{-5} \, \frac{\text{m}^d \, \text{s}}{\text{kg}}$ & $\text{mm}$ & $10^{-5}\,\text{m}$ & $10^{-3}\,\text{rad}$ & $10^{-5}\,\frac{\text{kg}\,\text{m}}{\text{s}^2}$ \\[0.2cm]
$2^{0}\phantom{^{-}}$ & $2^{8}\phantom{^{-}}$ & $16\phantom{.}$ & $1.550$ & $2.409$ & $9.372$ \\[0.2cm]
$2^{-1}$              & $2^{6}\phantom{^{-}}$ & $8\phantom{.5}$ & $6.811$ & $10.54$ & $41.01$ \\[0.2cm]
$2^{-2}$              & $2^{4}\phantom{^{-}}$ & $4\phantom{.5}$ & $23.74$ & $36.02$ & $142.5$ \\[0.2cm]
$2^{-3}$              & $2^{2}\phantom{^{-}}$ & $2\phantom{.5}$ & $56.85$ & $81.76$ & $340.3$ \\[0.2cm]
$2^{-4}$              & $2^{0}\phantom{^{-}}$ & $1\phantom{.5}$ & $96.80$ & $127.5$ & $577.0$ \\[0.2cm]
$2^{-5}$              & $2^{-2}$ & $0.5$ & $140.5$ & $171.2$ & $832.5$ \\[0.1cm]
\end{tabular}
\end{table}

\subsubsection{Test case 1C: mobility-induced divergence}
\label{sec:1C}

We have so far considered cases in which the mobility takes on a constant value independent of~$\varepsilon$, corresponding to a specific intrinsic slip length \cite{Jacqmin:2000kx}. However, it is known~\cite{Abels:2018ly,Jacqmin:2000kx,Demont_Stoter_van_Brummelen_2023} that to have a well-posed classical sharp interface limit, the mobility~$m$ must approach~$0$ 
alongside~$\varepsilon$, with scaling~$m \propto \varepsilon^\alpha$ with $0<\alpha<3$. To examine the effect of such a scaling of mobility, we will consider here the no-slip case in which the mobility scales quadratically with the interface thickness. The choice of $\alpha=2$ is substantiated by the fact that it provides a linear scaling of the corresponding diffusive slip length~$s_m$ as $\varepsilon \rightarrow +0$, although recent results indicate that a different scaling may be required to achieve optimal convergence to the sharp-interface 
limit~\cite{Demont_Stoter_van_Brummelen_2023}. We set a base parameter dependence according to $m=2.5 \!\times\! 10^{-6} \, \text{kg}^{-1} \, \text{m}^{d} \, \text{s}$, and thus a diffusive slip length $s_m = 0.5 \, \text{mm}$, for the smallest sampled interface thickness $\varepsilon = 5 \!\times\! 10^{-5} \, \text{m}$, along with all other parameter values depicted under test case 1C in Table~\TAB{model_param_vals}. Since the sharp interface limit without slip does not exist, we expect our results to diverge. We are interested in the rate of divergence, and the qualitative behavior of the (diverging) shear force as $\varepsilon, m \rightarrow +0$.

Table~\TAB{1C} presents the quantities of interest obtained from the aforementioned simulations. As expected, scaling the mobility according to $m \propto \varepsilon^2$, and hence the diffusive slip length per $s_m \propto \varepsilon$, results in increasing $\Delta x$, $\tilde{\vartheta}$, and $F_\textsc{S}$ as $\varepsilon,m \rightarrow +0$. The results indicate diverge in the limit, which is consistent with the classical contradiction for sharp-interface models subject to no-slip conditions, as presented in~\cite{huh1971hydrodynamic}. Moreover, the three quantities of interest diverge at a similar rate.

\begin{figure}
\centering
\includegraphics[width=0.7\textwidth]{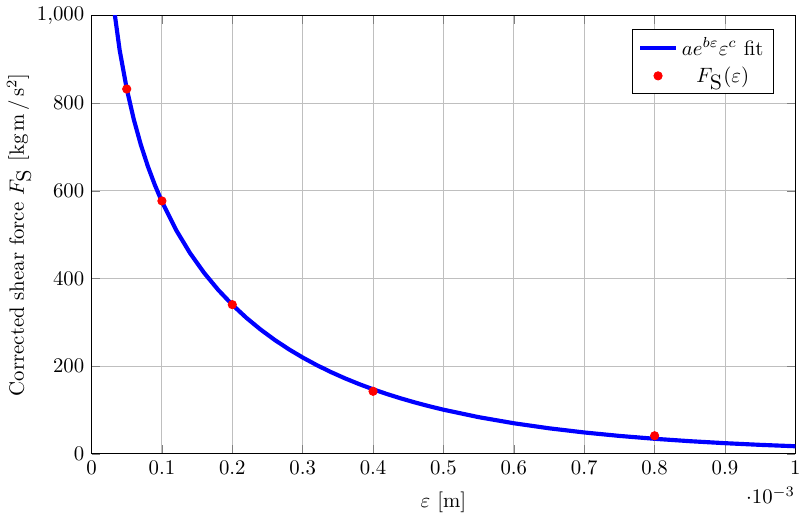}
\caption{Test case 1C: best fit for $ae^{b\varepsilon}\varepsilon^c$ to $F_\textsc{S}$ for the sampled $\varepsilon$ in the asymptotic regime with the three free parameters $a$, $b$, and $c$.}
\label{fig:1C_pol_fit}
\end{figure}

To elucidate the divergence of the results in Table~\TAB{1C}, Figure~\FIG{1C_pol_fit} present a fit to the computed shear-force data by means of an asymptotic exponential function, $ae^{b\varepsilon}\varepsilon^c$, . This fitted function is characterized by its two asymptotes for $\varepsilon \rightarrow \infty$ and $\varepsilon \rightarrow +0$, for which in the latter case the (reciprocal) factor becomes dominant. As the shear force $F_\textsc{S}$ has already been corrected for the Couette-flow-induced baseline contribution, there is no need for a constant parameter in the fitted model. As the fitting function, with $c \approx -0.3$, aligns closely with the computed shear force samples, it provides further support for the observed divergence as $\varepsilon \rightarrow +0$.


\subsection{Investigation of the sharp-interface limit with GNBC slip}
\label{ssec:numexp_GNBC}

The previous section established that when the fluid-fluid interface intersects a solid boundary at a contact point (2D) or line (3D), relying solely on the mobility parameter~$m$ to account for slip length leads to results that are inconsistent with the underlying sharp-interface model. On the one hand, if the diffuse slip length~$s_m$, and thus the mobility parameter $m$, is chosen as a fixed constant as $\varepsilon \rightarrow +0$, then the model does not converge to the classical sharp-interface model, yielding only approximately congruent interface dynamics; both near and away from the fluid--fluid interface and the contact points, the flow profiles differ drastically both qualitatively and quantitatively. On the other hand, when the mobility scales suitably~\cite{Abels:2018ly,Jacqmin:2000kx,Demont_Stoter_van_Brummelen_2023} with the interface thickness $\varepsilon$ as $\varepsilon, m \rightarrow +0$, no solution exists in the limit, as the diffusive slip length $s_m$ approaches zero as well. Slip must hence be accounted for in a different manner, to allow convergence to the classical sharp-interface limit with slip.

In this section, we consider slip via the Generalized Navier boundary condition, presented in Section~\ref{sec:wetting_BC}, and per~\cite{Stoter2023b}. In correspondence with Section~\ref{sec:1C}, we let $m \propto \varepsilon^2$ such that the NSCH sharp-interface limit converges to the classical sharp-interface model. However, in contrast to the previous section, we choose the baseline mobility $m$ small enough such that the diffusive slip length~$s_m$ is several orders of magnitude smaller than the GNBC slip length $s_\nu$, so that the former is negligible, to minimize the effect of the mobility-induced slip. We compare the results to the sharp-interface model subjected to the sharp-interface Generalized Navier boundary conditions~\eqref{eq:SI-GNBC} with the same slip length~$s_\nu$, and we investigate the convergence to the sharp-interface limit as~$\varepsilon, m \rightarrow +0$ for several slip lengths.

\subsubsection{Test case 2A: GNBC convergence analysis}
We investigate the convergence of the diffuse-interface model as $\varepsilon,m \rightarrow +0$ to the equivalent sharp-interface limit model for the case in which a finite slip is present through the generalized Navier boundary conditions on the wetting boundaries. Analogous to test case 1C in Section~\ref{sec:1C}, but with a smaller constant, we relate the mobility to the interface thickness according to $m\left(\varepsilon\right) = 0.4\,\varepsilon^2\,\textrm{kg}^{-1}\,\textrm{m}^{d-2}\,\textrm{s}$. With the remaining model parameters listed in Table~\TAB{model_param_vals} under test case~2A, the resulting diffusive slip length is~$s_m = 0.2\, \varepsilon$. We present a quantitative comparison of the contact-point displacement, macroscopic contact angle, and boundary shear force of the NSCH model with those of the corresponding sharp-interface model. Furthermore, a qualitative comparison is made of the flow fields of both models. As in the previous section, two separate GNBC slip lengths are regarded, viz.~$s_\nu = 2\,\textrm{mm}$ and~$s_\nu = 1\,\textrm{mm}$. For the considered parameter values, the diffusive slip length~$s_m$ for the largest examined~$\varepsilon$ is still approximately 3 times smaller than the smallest GNBC length-scale of $s_\nu = 1 \, \text{mm}$.

Table~\TAB{2A_comb} lists the computed quantities of interest for the aforementioned test case. The results in Table~\TAB{2A_comb} indicate that as~$\varepsilon$ decreases, each quantity of interest converges toward a limit point. By extrapolating the computed data and evaluating the result at $\varepsilon=0\,\textrm{m}$, a comparison can be made with the results of the sharp-interface model. The comparison yields good agreement between both models, with a relative discrepancy less than~$10^{-3}$ for all three quantities of interest and for both considered slip lengths, providing numerical corroboration of the the analysis in Section~\ref{sec:SI_lim_analysis}.
\begin{table}
\centering
\caption{Test case 2A: 
contact point displacements $\Delta x$, macroscopic contact angle $\tilde{\vartheta}$, and shear force~$F_\textsc{S}$ for the diffuse-interface model with GNBC and slip length $s_{\nu} = \{2,1\}\,\textrm{mm}$, for a decreasing sequence of interface-thickness parameters~$\varepsilon$ mobility~$m \propto\varepsilon^2$. The diffuse-interface results are extrapolated to $\varepsilon, m \rightarrow +0$ using a quadratic polynomial fit on the last three data points. Bottom row displays corresponding sharp-interface results.
\label{tab:2A_comb}}
\begin{tabular}{ccc|ccc|ccc}
 & & & \multicolumn{3}{c|}{$s_{\nu}=2\,\text{mm}$} & \multicolumn{3}{c}{$s_{\nu}=1\,\text{mm}$} \\[0.2cm]
\hline & & & & & & & & \\[-0.3cm]
$\varepsilon$ & $m$ & $s_m$ & $\Delta x$ & $\tilde{\vartheta}$ & $F_\textsc{S}$ & $\Delta x$ & $\tilde{\vartheta}$ & $F_\textsc{S}$ \\[0.1cm]
\hline
\hline & & & & & & & & \\[-0.3cm]
$10^{-4}\,\text{m}$ & $10^{-9} \, \frac{\text{m}^d \, \text{s}}{\text{kg}}$ & \textmu$\text{m}$ & $10^{-4}\,\text{m}$ & $10^{-2}\,\text{rad}$ & $10^{-3}\,\frac{\text{kg}\,\text{m}}{\text{s}^2}$ & $10^{-4}\,\text{m}$ & $10^{-2}\,\text{rad}$ & $10^{-3}\,\frac{\text{kg}\,\text{m}}{\text{s}^2}$ \\[0.2cm]
$2^{4}\phantom{^{-}}$ & $2^{10}$ & $320$ & $6.443$ & $8.255$ & $3.258$ & $8.957$ & $11.27$ & $4.924$ \\[0.2cm]
$2^{3}\phantom{^{-}}$ & $2^8\phantom{^{0}}$ & $160$ & $6.393$ & $8.074$ & $3.205$ & $9.140$ & $11.27$ & $4.984$ \\[0.2cm]
$2^{2}\phantom{^{-}}$ & $2^6\phantom{^{0}}$ & $\phantom{1}80$ & $6.291$ & $7.948$ & $3.143$ & $9.065$ & $11.15$ & $4.927$ \\[0.2cm]
$2^{1}\phantom{^{-}}$ & $2^4\phantom{^{0}}$ & $\phantom{1}40$ & $6.224$ & $7.882$ & $3.105$ & $8.962$ & $11.04$ & $4.866$ \\[0.2cm]
$2^{0}\phantom{^{-}}$ & $2^2\phantom{^{0}}$ & $\phantom{1}20$ & $6.192$ & $7.852$ & $3.088$ & $8.898$ & $10.98$ & $4.829$ \\[0.2cm]
$2^{-1}$              & $2^0\phantom{^{0}}$ & $\phantom{1}10$ & $6.178$ & $7.840$ & $3.081$ & $8.870$ & $10.96$ & $4.813$ \\[0.2cm]
\hdashline & & & & & & & & \\[-0.3cm]
\multicolumn{3}{c|}{extrapolated ($\varepsilon \rightarrow +0$)} & $6.167$ & $7.831$ & $3.075$ & $8.845$ & $10.93$ & $4.800$ \\[0.2cm]
\hline & & & & & & & & \\[-0.3cm]
\multicolumn{3}{c|}{NS w/slip}             & $6.171$ & $7.836$ & $3.077$ & $8.848$ & $10.94$ & $4.801$ \\[0.1cm]
\end{tabular}
\end{table}

\begin{figure}
\centering
\includegraphics[width=1.00\textwidth]{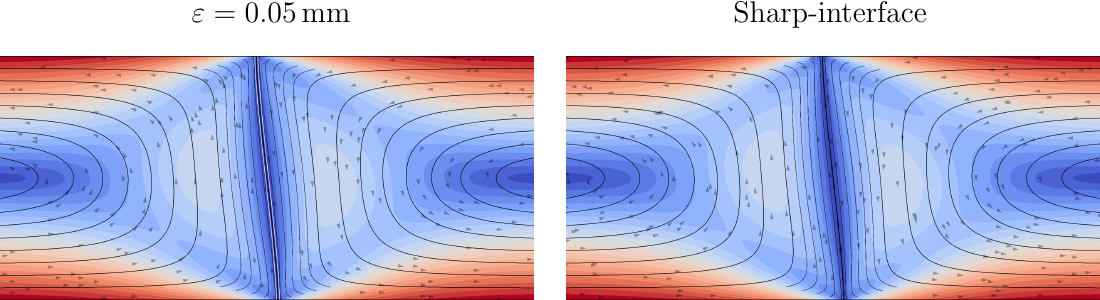}
\caption{Velocity fields and stream lines for the NSCH model with diffuse-interface GNBC~\eqref{eq:wetting_BC_sub2} ({\em left\/}) and the sharp-interface model with GNBC~\eqref{eq:SI-GNBC} ({\em right\/}) for slip length~$s_\nu = 1\,\textrm{mm}$.
\label{fig:2A_plot}}
\end{figure}

Figure~\FIG{2A_plot} juxtaposes the computed flow fields for the diffuse-interface model with GNBC~\eqref{eq:wetting_BC_sub2}, and the sharp-interface model with GNBC~\eqref{eq:SI-GNBC}, for slip length~$s_\nu = 1\,\textrm{mm}$. It is to be noted that the right plot in Figure~\FIG{2A_plot} coincides with the bottom-right plot in Figure~\ref{fig:1A_velocity_fields}, and has been repeated here to facilitate a direct comparison. For the diffuse-interface model, the flux through the interface and the fluid circulation near the center of and through the interface, which were observed in the diffuse-interface model with no-slip boundary conditions and mobility-mediated slip in Section~\ref{sec:1A}, have essentially disappeared. The velocity fields for the diffuse-interface model and the sharp-interface model in Figure~\FIG{2A_plot} are indistignuishable, except in the vicinity of the  fluid--fluid interface, where the finite interface thickness is visible. The results indicate that if the diffuse-interface model is equipped with the GNBC with finite slip length~$s_\nu > 0$, and the mobility decreases as $m \propto \varepsilon^\alpha$ with $0 < \alpha < 3$ as $\varepsilon\to{}+0$, the diffuse-interface solution converges to the corresponding sharp-interface solution. It is noteworthy that the close alignment between the diffuse-interface and sharp-interface results, which have been obtained from different code-bases, suggests that the presented test case and corresponding results can be regarded as a benchmark for validating binary-fluid flow models with wetting.


\subsection{Triple-wedge flow}
\label{sssec:2BC}
Finally, we seek to reproduce the triple-wedge flows as documented in~\cite{huh1971hydrodynamic} by means of the diffuse-interface model with GNBC. In Ref.~\cite{huh1971hydrodynamic}, it is shown there exists a critical viscosity ratio $R_\textsc{c}$, dependent on the microscopic contact angle $\vartheta$, for which the radial velocity component (with the fluid--fluid--solid contact point defined as the origin) on the fluid--fluid interface in a region around the contact point is exactly zero. This critical viscosity ratio is given by
\begin{equation}
R_\textsc{c}(\vartheta) \coloneqq \frac{(\vartheta\cos(\vartheta)-\sin(\vartheta))((\vartheta-\pi)^2-(\sin(\vartheta))^2)}{((\vartheta-\pi)\cos(\vartheta)-\sin(\vartheta))((\vartheta)^2-(\sin(\vartheta))^2)}\,.
\end{equation}
Throughout all preceding test cases, a $\pi/2\,\textrm{rad}$ microscopic contact angle has been imposed, for which the corresponding critical viscosity ratio is exactly $R_\textsc{c}(\pi/2) = 1$. Since in the preceding simulations the fluid properties have been fully symmetric, the fluid--fluid viscosity ratios are exactly equal to the critical viscosity ratio $R_\textsc{c}$, and one can indeed observe that throughout all the relevant velocity field plots, the radial velocity component along the fluid--fluid interface has been approximately 0. However, when the viscosity ratio $\eta_\AA / \eta_\LL$ drops below the critical viscosity ratio $R_\textsc{c}$, then the radial velocity component along the fluid--fluid interface is directed toward the upper contact point~\cite{huh1971hydrodynamic}, invoking a triple-wedge-flow configuration. The middle wedge develops on the side of the ambient phase $\AA$, where $\varphi = -1$. Conversely, when the viscosity ratio $\eta_\AA / \eta_\LL$ is larger than $R_\textsc{c}$, the velocity direction is reversed, pointing downward toward the lower contact point, and the middle wedge is located on the side of the liquid phase $\LL$ ($\varphi = 1$). In the remainder of this section, we prompt the triple-wedge flow in two ways: by changing the ambient viscosity $\eta_\AA$ such that the fluid--fluid viscosity ratio subseeds the critical viscosity ratio, and by decreasing the microscopic contact angle $\vartheta$ such that the critical viscosity ratio~$R_\textsc{c}$ exceeds the unitary fluid--fluid viscosity ratio.

\begin{figure}
\centering
\includegraphics[width=1.00\textwidth]{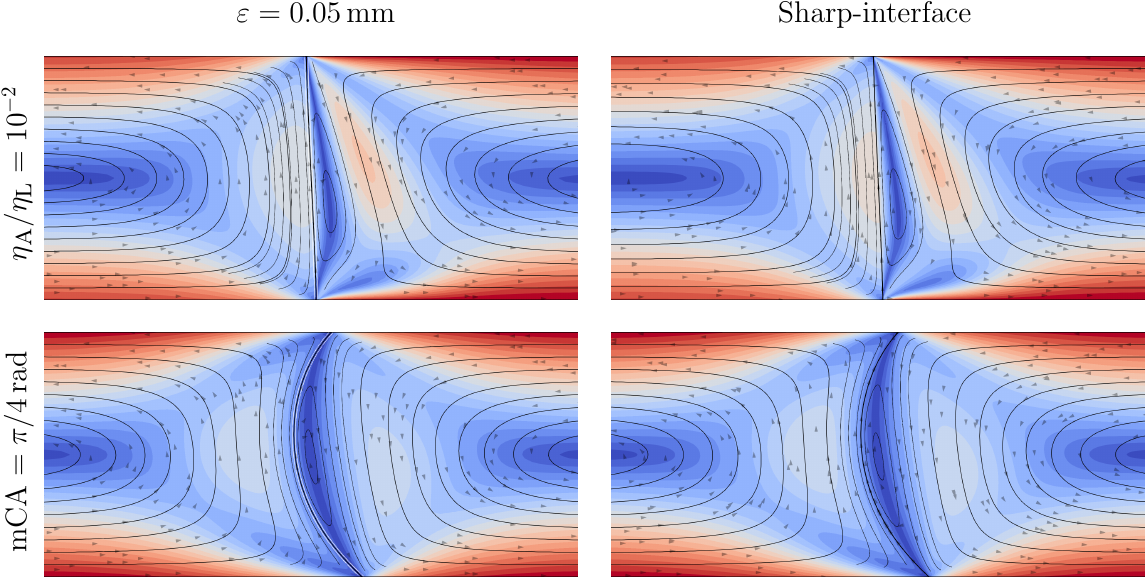}
\caption{Velocity fields and streamlines for the NSCH (left) and sharp-interface (right) models with viscosity ratio $\eta_\AA / \eta_\LL = 10^{-2}$ (top) and microscopic contact angle $\vartheta = \pi/4\,\textrm{rad}$ (bottom).
\label{fig:2B_2C_triple_wedge_visc_45deg}}
\end{figure}

Figure~\FIG{2B_2C_triple_wedge_visc_45deg} ({\em top\/}) presents the flow fields corresponding to the numerical experiment with parameter values as listed under test case 2B in Table~\TAB{model_param_vals}. For the figures in the top row, the viscosity of the ambient fluid has been reduced by a factor~$100$, representing water in a regular room-temperature environment. The viscosity ratio $\eta_\AA / \eta_\LL = 10^{-2}$ is much smaller than the critical viscosity ratio $R_\textsc{c}(\pi/2) = 1$, generating a triple-wedge flow with the middle wedge located in the ambient phase~$\AA$. One can observe a clear asymmetry in the flow fields. It is to be noted that due to the decrease of the viscosity in the ambient fluid phase, there is a decrease of the slip length~$s_\nu$ (as well as $s_m$) in the ambient phase to $s_\nu = 10\,\textrm{\textmu}m$, dropping below the interface thickness parameter $\varepsilon = 50\,\textrm{\textmu}m$. The effects of the small slip length are visible close to the contact points, where relatively high velocity gradients occur, and away from the contact points where there is a resulting higher absolute velocity. Even though the finite interface thickness interferes with the resulting small slip length in the ambient phase, the triple wedge flow --- in a region close to the contact points, but away from the diffuse-interface --- is clearly observable. Despite the complexity of the fluid dynamics in this test case, the diffuse-interface and sharp-interface results still display a close resemblance.

Finally, Figure~\FIG{2B_2C_triple_wedge_visc_45deg} ({\em bottom\/}) displays the flow fields for test case 2C in Table~\TAB{model_param_vals}. We once again impose a unitary fluid--fluid viscosity ratio $\eta_\AA / \eta_\LL = 1$, and invoke the triple-wedge flow by means of increasing the critical viscosity ratio $R_\textsc{c}$. By decreasing the microscopic contact angle to half its original value, $\vartheta = \pi/4\,\textrm{rad}$, we effectively increase the critical viscosity ratio to $R_\textsc{c}(\pi/4) \approx 2.76433$, ensuring the actual viscosity ratio is below the critical viscosity ratio. As a result, we observe the three wedges to have the same flow directions and be in the same fluid phases as is the case for test case 2B, albeit with drastically different flow profiles away from the contact points. The diffuse-interface and sharp-interface results again display close agreement for this complex case.

The results of test cases 2B and 2C convey that the diffuse-interface model equipped with the diffuse-interface GNBC~\eqref{eq:wetting_BC_sub2} allows to capture complex wetting behavior, in agreement with the underlying sharp-interface model with sharp-interface GNBC~\eqref{eq:SI-GNBC}.


\section{Conclusions}
\label{sec:concl}
The implicit representation of the interface in diffuse-interface binary-fluid models holds considerable promise in describing complex phenomena within fluid mechanics, including topological changes of the fluid–fluid interface and dynamic wetting. In the absence of topological changes in the interface, diffuse-interface models are expected to converge to their corresponding classical sharp-interface counterparts, known as the sharp-interface limit. This convergence occurs as the interface-thickness parameter, $\varepsilon$, and the mobility parameter, $m \propto \varepsilon^\alpha$, pass to zero, where a proper relative scaling with $0 < \alpha < 3$ is required. Contemporary understanding of the sharp-interface limit in the case of dynamic wetting, i.e. if the fluid--fluid interface intersects an adjacent solid surface and the corresponding contact line exhibits motion relative to that surface, is incomplete. In this article, we investigated the limit behavior of the Abels--Garcke--Gr\"{u}n NSCH model in dynamic-wetting scenarios. In cases of dynamic wetting, sharp-interface binary-fluid models subject to a no-slip condition exhibit a non-integrable stress singularity at the contact line. To avoid such a singularity, a slip mechanism must be introduced, the standard option being a generalized Navier boundary condition. 

Diffuse-interface binary-fluid models possess an intrinsic slip mechanism, which has been propounded as a means to introduce slip into diffuse-interface binary-fluid models equipped with no-slip boundary conditions. However, the intrinsic slip length~$s_{m}$ is proportional to~$\sqrt{m}$ and, hence, a finite slip length precludes the limit $m\to{}+0$. By means of numerical experiments for a Couette test case, for the diffuse-interface model and for a corresponding sharp-interface model with a generalized Navier boundary condition with slip length~$s_{\nu}$, we established that if $s_m=s_{\nu}$, the diffuse-interface model and sharp-interface model exhibit very similar interface profiles. However, in passing the interface thickness in the diffuse-interface model to zero, the flow field in the diffuse-interface model is distinctly different from the flow field in the sharp-interface model. This is consistent with the fact that for fixed mobility~$m$, in the limit $\varepsilon\to{}+0$, the diffuse-interface model does not converge to the classical sharp-interface model, but to the Navier--Stokes/Mullins--Sekerka model, instead. If, in contrast, $m\propto\varepsilon^2$ such that $m\to{}+0$ as $\varepsilon\to+0$ to recover the classical sharp-interface binary-fluid model then, indeed, the shear force, contact-point displacement and macroscopic contact angle in the diffuse-interface model diverge in the limit, consistent with the degeneration of the underlying sharp-interface model with no-slip conditions. The mobility-mediated slip in diffuse-interface models is therefore severely restricted in modeling dynamic-wetting phenomena.

To resolve the aforementioned issues associated with no-slip boundary conditions in dynamic-wetting scenarios, we equipped the NSCH model with a generalized Navier boundary condition and a dynamic contact-angle condition. We derived that the generalized Navier boundary condition and dynamic (or static) contact-angle conditions in the diffuse-interface model are thermodynamically admissible, in the sense that these provide non-negative free-energy dissipation. By means of a formal analysis, we established that the generalized Navier boundary condition for the diffuse-interface model converges to its sharp-interface counterpart as the interface thickness passes to zero. By virtue of the slip provided by the generalized Navier boundary condition, the mobility in the diffuse-interface model can be selected independent of the slip length, allowing~$m$ to vanish in tandem with~$\varepsilon\to+0$.
By means of systematic numerical experiments for a Couette test case, including comparisons to corresponding sharp-interface simulations, we established that the solution of the diffuse-interface model equipped with the generalized Navier boundary condition converges to the classical sharp-interface solution as~$\varepsilon, m \rightarrow +0$. As a byproduct of the cross-validation of the diffuse-interface and sharp-interface results, obtained by means of different code bases, the presented results can serve as a benchmark for future investigations of dynamic-wetting simulations.
Finally, we have confirmed that the diffuse-interface model with generalized Navier boundary condition is capable of reproducing the triple-wedge flow patterns that are known to arise in a region around the contact points when the fluid--fluid viscosity ratio deviates from the critical viscosity ratio --- either by altering the viscosity ratio, or changing the microscopic contact angle.


\section*{Acknowledgments}
This research was partly conducted within the Industrial Partnership Program {\it Fundamental Fluid Dynamics Challenges in 
Inkjet Printing\/} ({\it FIP\/}), a joint research program of Canon Production Printing, Eindhoven University of Technology,
University of Twente, and the Netherlands Organization for Scientific Research (NWO). T.H.B.\ Demont gratefully acknowledges financial support through the FIP program.

\bibliographystyle{plain}
\bibliography{wetting}

\begin{thebibliography}{10}

\bibitem{Abels:2018ly}
H.~Abels and H.~Garcke.
\newblock {\em Weak Solutions and Diffuse Interface Models for Incompressible Two-Phase Flows}, pages 1267--1327.
\newblock Springer International Publishing, Cham, 2018.

\bibitem{Abels:2012vn}
H.~Abels, H.~Garcke, and G.~Gr{\"u}n.
\newblock Thermodynamically consistent, frame indifferent diffuse interface models for incompressible two-phase flows with different densities.
\newblock {\em Math. Mod. Meth. Appl. Sci.}, 22:1150013, 2012.

\bibitem{Abels:2014ca}
H.~Abels and D.~Lengeler.
\newblock On sharp interface limits for diffuse interface models for two-phase flows.
\newblock {\em Interfaces Free Boundaries}, 395--418, 2014.

\bibitem{Arrhenius:1887xr}
S.~Arrhenius.
\newblock {\"U}ber die innere {R}eibung verd{\"u}nnter w{\"a}sseriger {L}{\"o}sungen.
\newblock {\em Z. Physik. Chem.}, 1U(1):285--298, 1887.

\bibitem{Barrett:1999nx}
J.~Barrett, J.~Blowey, and H.~Garcke.
\newblock Finite element approximation of the {C}ahn--{H}illiard equation with degenerate mobility.
\newblock {\em SIAM J. Numer. Anal.}, 37(1):286--318, 2016/11/21 1999.

\bibitem{Bonart:2019re}
H.~Bonart, C.~Kahle, and J.-U. Repke.
\newblock Comparison of energy stable simulation of moving contact line problems using a thermodynamically consistent {C}ahn--{H}illiard {N}avier--{S}tokes model.
\newblock {\em J. Comput. Phys.}, 399:108959, 2019.

\bibitem{Cox:1986jq}
R.G. Cox.
\newblock The dynamics of the spreading of liquids on a solid surface. part 1. viscous flow.
\newblock {\em Journal of Fluid Mechanics}, 168:169--194, 1986.

\bibitem{Demont_Stoter_van_Brummelen_2023}
T.H.B. Demont, S.K.F. Stoter, and E.H. van Brummelen.
\newblock Numerical investigation of the sharp-interface limit of the {N}avier--{S}tokes--{C}ahn--{H}illiard equations.
\newblock {\em Journal of Fluid Mechanics}, 970:A24, 2023.

\bibitem{Demont:2022dk}
T.H.B. Demont, G.J. van Zwieten, C.~Diddens, and E.H. van Brummelen.
\newblock A robust and accurate adaptive approximation method for a diffuse-interface model of binary-fluid flows.
\newblock {\em Comput. Methods Appl. Mech. Engrg.}, 400:115563, 2022.

\bibitem{diddens2023bifurcation}
C.~Diddens and D.~Rocha.
\newblock Bifurcation tracking on moving meshes and with consideration of azimuthal symmetry breaking instabilities, 2023.

\bibitem{Gerbeau:2009km}
J.{-}F. Gerbeau and T.~Leli{\`e}vre.
\newblock Generalized {N}avier boundary condition and geometric conservation law for surface tension.
\newblock {\em Computer Methods in Applied Mechanics and Engineering}, 198(5):644--656, 2009.

\bibitem{Grun:2016gi}
G.~Gr\"un, F.~Guill\'en-Gonz\'alez, and S.~Metzger.
\newblock On fully decoupled, convergent schemes for diffuse interface models for two-phase flow with general mass densities.
\newblock {\em Commun. Comput. Phys.}, 19(5):1473--1502, 2016.

\bibitem{Hohenberg:1977hh}
P.C. Hohenberg and B.I. Halperin.
\newblock Theory of dynamic critical phenomena.
\newblock {\em Rev. Mod. Phys.}, 49:435--479, 1977.

\bibitem{huh1971hydrodynamic}
C.~Huh and L.~E. Scriven.
\newblock Hydrodynamic model of steady movement of a solid-liquid-fluid contact line.
\newblock {\em Journal of Colloid and Interface Science}, 35:85--101, 1971.

\bibitem{Jacqmin:2000kx}
D.~Jacqmin.
\newblock Contact-line dynamics of a diffuse fluid interface.
\newblock {\em J. Fluid Mech.}, 402:57--88, 2000.

\bibitem{Khatavkar:2006gf}
V.V. Khatavkar, P.D. Anderson, and H.E.H. Meijer.
\newblock On scaling of diffuse--interface models.
\newblock {\em Chem. Eng. Sci.}, 61(8):2364--2378, 2006.

\bibitem{lohse2022fundamental}
D.~Lohse.
\newblock Fundamental fluid dynamics challenges in inkjet printing.
\newblock {\em Annual review of fluid mechanics}, 54(1):349--382, 2022.

\bibitem{Lowengrub:1998uq}
J.~Lowengrub and L.~Truskinovsky.
\newblock Quasi-incompressible {C}ahn-{H}illiard fluids and topological transitions.
\newblock {\em Proceedings of the Royal Society of London A: Mathematical, Physical and Engineering Sciences}, 454:2617--2654, 1998.

\bibitem{rougier2021slip}
V.~Rougier, J.~Cellier, M.~Gomina, and J.~Br{\'e}ard.
\newblock Slip transition in dynamic wetting for a generalized {N}avier boundary condition.
\newblock {\em Journal of Colloid and Interface Science}, 583:448--458, 2021.

\bibitem{Seppecher:1996kx}
P.~Seppecher.
\newblock Moving contact lines in the {C}ahn-{H}illiard theory.
\newblock {\em Int. J. Engng. Sci.}, 34:977--992, 7 1996.

\bibitem{Simsek:2018gb}
M.~Shokrpour~Roudbari, G.~\c{S}im\c{s}ek, E.H. {van Brummelen}, and K.G. {van der Zee}.
\newblock Diffuse-interface two-phase flow models with different densities: A new quasi-incompressible form and a linear energy-stable method.
\newblock {\em Math. Models Methods Appl. Sci.}, 28:733--770, 2018.

\bibitem{Shokrpour-Roudbari:2016dp}
M.~Shokrpour~Roudbari, E.~H. {van Brummelen}, and C.V. Verhoosel.
\newblock A multiscale diffuse-interface model for two-phase flow in porous media.
\newblock {\em Computers \& Fluids}, 141:212--222, 2016.

\bibitem{Stoter2023b}
S.K.F. Stoter, T.B. van Sluijs, T.H.B. Demont, E.H. van Brummelen, and C.V. Verhoosel.
\newblock {Stabilized immersed isogeometric analysis for the Navier–Stokes–Cahn–Hilliard equations, with applications to binary-fluid flow through porous media}.
\newblock {\em Computer Methods in Applied Mechanics and Engineering}, 417B:116483, 2023.

\bibitem{Brummelen:2021aw}
E.H. van Brummelen, T.H.B. Demont, and G.J. van Zwieten.
\newblock An adaptive isogeometric analysis approach to elasto-capillary fluid-solid interaction.
\newblock {\em Int. J. Numer. Meth. Engng.}, 122(19):5331--5352, 2021/09/30 2021.

\bibitem{Brummelen:2016qa}
{E.H.} {van Brummelen}, H.~Shokrpour~Roudbari, and {G.J.} {van Zwieten}.
\newblock Elasto-capillarity simulations based on the {N}avier-{S}tokes-{C}ahn-{H}illiard equations.
\newblock In {\em Advances in Computational Fluid-Structure Interaction and Flow Simulation}, Modeling and Simulation in Science, Engineering and Technology, pages 451--462. Birkh{{\"a}}user, 2016.

\bibitem{Yue:2011uq}
P.~Yue and J.J. Feng.
\newblock Wall energy relaxation in the {C}ahn--{H}illiard model for moving contact lines.
\newblock {\em Phys. Fluids}, 23:012106, 2011.

\bibitem{Yue:2010hq}
P.~Yue, C.~Zhou, and J.J. Feng.
\newblock Sharp-interface limit of the {C}ahn--{H}illiard model for moving contact lines.
\newblock {\em J. Fluid Mech.}, 645:279--294, 2010.

\end{thebibliography}

\end{document}